\ifdefined\synctex\synctex=1\relax\fi
\documentclass[preprint,11pt,a4paper]{elsarticle}

\usepackage{amsmath,amssymb,amsthm,mathtools,bm}
\numberwithin{equation}{section}
\usepackage{geometry}
\usepackage{enumitem}
\usepackage{xcolor}
\usepackage[hidelinks]{hyperref}
\allowdisplaybreaks
\newtheorem{theorem}{Theorem}[section]
\newtheorem{lemma}[theorem]{Lemma}
\newtheorem{proposition}[theorem]{Proposition}
\newtheorem{corollary}[theorem]{Corollary}
\theoremstyle{remark}
\newtheorem{remark}[theorem]{Remark}

\newcommand{\Om}{\Omega}
\newcommand{\R}{\mathbb R}
\newcommand{\dd}{\,\mathrm d}
\newcommand{\Tmax}{T_{\max}}
\newcommand{\Bunit}{|B_1(0)|}
\newcommand{\uz}{\underline z}
\newcommand{\F}{\mathcal F}

\begin{document}

\begin{frontmatter}

\title{Global Boundedness and Finite-Time Blow-up for a
Quasilinear Chemotaxis--May--Nowak System}

\author[nwpu]{Xueli Bai}
\ead{xlbai2015@nwpu.edu.cn}
\author[sysu]{Fang Li}
\ead{lifang55@mail.sysu.edu.cn}
\author[nwpu]{Suying Liu\texorpdfstring{\corref{cor1}}{}}
\ead{suyingliu@nwpu.edu.cn}
\cortext[cor1]{Corresponding author.}

\address[nwpu]{School of Mathematics and Statistics, Northwestern
Polytechnical University, 1 Dongxiang Road, Chang'an District, Xi'an,
Shaanxi 710129, P.R. China}
\address[sysu]{School of Mathematics, Sun Yat-sen University, 135 Xingang
West Road, Haizhu District, Guangzhou 510275, P.R. China}

\begin{abstract}
This paper continues the work of Bai, Li, Liu and Luo \cite{BaiLiLiuLuoPreprint} on a
quasilinear chemotaxis--May--Nowak system. We establish global boundedness
in three dimensions for both the parabolic--elliptic--parabolic and fully
parabolic systems when $m>4/3$. For the former, we construct finite-time
blow-up solutions in balls when $n=3$ and $1<m\le4/3$, and when $n\ge4$ and
$0<m\le2-2/n$.
\end{abstract}

\end{frontmatter}

\section{Introduction and main results}

The classical May--Nowak model describes within-host viral dynamics
\cite{NowakBangham1996,NowakMay2000}. Stancevic et al.
\cite{Stancevic2013} incorporated diffusion and chemotaxis, and the
resulting framework was further developed in
\cite{BellomoPainterTaoWinkler2019,TaoWinkler2021}. Building on these
works, Wang and Wang \cite{wang2025} extended the model to nonlinear
diffusion. We consider the corresponding quasilinear system
in a bounded domain $\Om\subset\R^n$ with smooth boundary:
\begin{equation}\label{sys}
\left\{
\begin{aligned}
 u_t&=\nabla\!\cdot\!\bigl((1+u)^{m-1}\nabla u\bigr)
      -\chi\nabla\!\cdot(u\nabla v)-u-uw+\varphi,
      &&x\in\Om,\ t>0,\\
 \tau v_t&=\Delta v-v+uw,
      &&x\in\Om,\ t>0,\\
 w_t&=\Delta w-w+v,
      &&x\in\Om,\ t>0,\\
 \partial_\nu u&=\partial_\nu v=\partial_\nu w=0,
      &&x\in\partial\Om,\ t>0,\\
 u(x,0)&=u_0(x),\
 \tau v(x,0)=\tau v_0(x),\
 w(x,0)=w_0(x),
      &&x\in\Om.
\end{aligned}
\right.
\end{equation}

Here $u$, $v$, and $w$ denote the densities of uninfected cells, infected
cells, and free virus particles, respectively, and $\varphi$ is the cell
supply. The exponent $m$ describes density-dependent diffusion.
The choices $\tau=0$ and $\tau=1$ correspond to the
parabolic--elliptic--parabolic and fully parabolic systems, respectively.

For results on these systems and related models, see
\cite{BellomoPainterTaoWinkler2019,TaoWinkler2021,wang2025,Winkler2019Acta,ZhengWang2026}.
This paper continues our previous work with Luo \cite{BaiLiLiuLuoPreprint}. We determine the optimal boundedness--blow-up
threshold for the three-dimensional parabolic--elliptic--parabolic system,
establish boundedness in the same range for the fully parabolic system,
and prove finite-time blow-up in higher dimensions.

Here $\chi>0$, $m>0$, and
\begin{equation}\label{source-ass}
 \varphi\in C^\alpha(\overline\Om\times[0,\infty)),
 \qquad 0\le\varphi(x,t)\le\varphi_*
 \quad (x\in\Om,\ t>0)
\end{equation}
for some $\varphi_*>0$ and $\alpha\in(0,1)$.  The initial data satisfy
\begin{equation}\label{data-ass}
 \begin{gathered}
 u_0,w_0\in W^{2,\infty}(\Om),\qquad u_0,w_0\ge0,\\
 u_0\not\equiv0,\quad w_0\not\equiv0,
 \qquad \partial_\nu u_0=\partial_\nu w_0=0.
 \end{gathered}
\end{equation}
When $\tau=1$, we assume in addition that
\begin{equation}\label{data-full}
 v_0\in W^{2,\infty}(\Om),\qquad v_0\ge0,
 \qquad \partial_\nu v_0=0.
\end{equation}

For either value of $\tau$, standard local theory yields a nonnegative
classical solution on a maximal interval $[0,\Tmax)$ satisfying the
extensibility criterion
\begin{equation*} \Tmax<\infty
 \quad\Longrightarrow\quad
 \limsup_{t\nearrow\Tmax}
 \bigl(\|u(\cdot,t)\|_{L^\infty}
       +\|w(\cdot,t)\|_{L^\infty}\bigr)=\infty;
\end{equation*}
see, for example, \cite[Lemma~1.1]{wang2025}.

Our main results are as follows.

\begin{theorem}
\label{thm:main-elliptic-global}
Let $n=3$ and $\tau=0$, and suppose that \eqref{source-ass} holds.
If $m>4/3$, then, for every pair $(u_0,w_0)$ satisfying
\eqref{data-ass}, the corresponding solution of \eqref{sys} is global and
uniformly bounded; more precisely,
\[
 \sup_{t>0}\left(
 \|u(\cdot,t)\|_{L^\infty(\Om)}
 +\|v(\cdot,t)\|_{W^{1,\infty}(\Om)}
 +\|w(\cdot,t)\|_{L^\infty(\Om)}
 \right)<\infty.
\]
\end{theorem}

\begin{theorem}
\label{thm:main-elliptic-blowup}
Let $n=3$, $\tau=0$, and $\Om=B_R(0)\subset\R^3$, and assume in addition that $\varphi$ is radial in
space.  If $1<m\le4/3$, then, for every $\mu>0$, there exist nonnegative
radial initial data $(u_0,w_0)$ satisfying \eqref{data-ass} and
$\int_\Om u_0=\mu$ such that the corresponding solution of \eqref{sys}
blows up in finite time in the sense that
\[
 \Tmax<\infty,
 \qquad
 \limsup_{t\nearrow\Tmax}\|u(\cdot,t)\|_{L^\infty(\Om)}=\infty.
\]
More precisely, one may take $w_0\equiv W_0$ with arbitrary $W_0>0$ when
$1<m<4/3$, whereas at $m=4/3$ one may take any sufficiently large $W_0>0$.
\end{theorem}

\begin{theorem}
\label{thm:main-high-blowup}
Let $n\ge4$, $\Om=B_R(0)\subset\R^n$, $\tau=0$, and
$0<m\le 2-2/n$.  Suppose that \eqref{source-ass} holds and that
$\varphi$ is radial in space.  For every $\mu>0$, there exist
nonnegative radial initial data $(u_0,w_0)$ satisfying \eqref{data-ass},
with $\int_\Om u_0=\mu$ such that the
corresponding solution of \eqref{sys} satisfies
\[
 \Tmax<\infty,
 \qquad
 \limsup_{t\nearrow\Tmax}\|u(\cdot,t)\|_{L^\infty(\Om)}=\infty.
\]
One may take $w_0\equiv W$ with any $W>0$ if $m<2-2/n$
and sufficiently large $W>0$ if $m=2-2/n$.
\end{theorem}

\begin{theorem}
\label{thm:main-parabolic}
Let $n=3$, $\tau=1$ and $m>4/3$, and suppose that \eqref{source-ass},
\eqref{data-ass}, and \eqref{data-full} hold.  Then the corresponding
solution of \eqref{sys} is global and uniformly bounded; more precisely,
\[
 \sup_{t>0}\left(
 \|u(\cdot,t)\|_{L^\infty(\Om)}
 +\|v(\cdot,t)\|_{W^{1,\infty}(\Om)}
 +\|w(\cdot,t)\|_{L^\infty(\Om)}
 \right)<\infty.
\]
\end{theorem}

\begin{remark}
The arguments can be adapted to quasilinear versions of the general
reaction systems in \cite{TaoWinkler2021}. For global boundedness, the adaptation uses
additional dissipative upper bounds on the reactions and subquadratic
growth of the infection rate in $w$. The reaction bounds in that work
suffice for three-dimensional blow-up when $1<m<4/3$; at $m=4/3$, the
construction also uses a sufficiently strong lower bound on the infection rate. The
higher-dimensional construction further requires suitable control of
cell loss, infection rates comparable to $w$, and nondecreasing concave virus
production with convex virus loss.
\end{remark}

The remainder of the paper is organized as follows.
Section~\ref{sec:tools} collects preliminary results and analytic tools.
Sections~\ref{sec:elliptic} and~\ref{sec:blowup-subcritical} establish
global boundedness and finite-time blow-up via subsolution construction,
respectively, for $n=3$ and $\tau=0$, identifying $4/3$ as the optimal
critical exponent.
Section~\ref{sec:high-blowup} proves
Theorem~\ref{thm:main-high-blowup} on finite-time blow-up in higher
dimensions by introducing a new energy structure, and
Section~\ref{sec:full} proves Theorem~\ref{thm:main-parabolic} on global
boundedness for the three-dimensional fully parabolic system.

\section{Preliminaries}\label{sec:tools}

First of all, we present the known local solvability,
extensibility property and $L^1$ estimates used throughout the paper for the convenience of readers.
See \cite[Lemma 1.1]{wang2025} for example.

\begin{lemma}\label{global_bounded_lemma}
Let $\tau\in\{0,1\}$, $n\ge2$, and $m>0$. Assume \eqref{source-ass}, and
suppose that the initial data satisfy \eqref{data-ass} and, when $\tau=1$,
also \eqref{data-full}. Then there exist
$T_{\max}\in(0,\infty]$ and a nonnegative classical solution $(u,v,w)$ of
\eqref{sys} on $\Omega\times(0,T_{\max})$ such that
$$
 u,w\in C^0\bigl(\overline\Omega\times[0,T_{\max})\bigr)
 \cap C^{2,1}\bigl(\overline\Omega\times(0,T_{\max})\bigr),
$$
and
$$
 v\in
 \begin{cases}
 C^{2,0}\bigl(\overline\Omega\times(0,T_{\max})\bigr),&\tau=0,\\
 C^0\bigl(\overline\Omega\times[0,T_{\max})\bigr)
 \cap C^{2,1}\bigl(\overline\Omega\times(0,T_{\max})\bigr),&\tau=1.
 \end{cases}
$$
Moreover, $u,v,w>0$ in $\overline\Omega\times(0,T_{\max})$. If
$T_{\max}<\infty$, then
$$
 \limsup_{t\nearrow T_{\max}}
 \left(
 \|u(\cdot,t)\|_{L^\infty(\Omega)}
 +\|w(\cdot,t)\|_{L^\infty(\Omega)}
 \right)=\infty.
$$
\end{lemma}

For every $m>0$, the coefficient $(1+u)^{m-1}$ is smooth and
bounded away from zero on bounded intervals of $u$, so the local theory
also applies when $0<m\le1$.

We first collect the common mass estimate and the analytic tools used below.
For the mass estimate in Lemma~\ref{lem:mass}, see
\cite{BaiLiLiuLuoPreprint} and \cite[Lemma~1.1]{wang2025}.

\begin{lemma}\label{lem:mass}
Let $\tau\in\{0,1\}$, and suppose that $(u,v,w)$ is a solution of
\eqref{sys}.  Then there exists $C>0$ such that
\begin{equation*} \sup_{t\in(0,\Tmax)}
 \left(
  \|u(\cdot,t)\|_{L^1(\Om)}
  +\|w(\cdot,t)\|_{L^1(\Om)}
 \right)\le C.
\end{equation*}
\end{lemma}

\begin{proof}
Integrating the equations in \eqref{sys} over $\Om$ and suitably adding
the resulting identities, we obtain the assertion by standard ODE comparison.
\end{proof}

\begin{lemma}\label{lem:elliptic}
Let $n=3$. If $-\Delta z+z=f$ in $\Om$ and $\partial_\nu z=0$, then, for every
$q\in(1,3)$,
\begin{equation*} \|z\|_{L^q(\Om)}\le C_q\|f\|_{L^1(\Om)}.
\end{equation*}
\end{lemma}

Lemma~\ref{lem:elliptic} follows directly from the
Brezis--Strauss estimate \cite{BrezisStrauss1973}.

We also record the two damped Neumann heat-semigroup estimates used below;
see \cite[Lemma~1.3]{Winkler2010JDE}.

\begin{lemma}
\label{lem:neumann-semigroup}
Let $n=3$ and let $\Delta_N$ denote the realization of the Laplacian in $\Om$ subject to
homogeneous Neumann boundary conditions.  Then
\begin{equation*} \|e^{t(\Delta_N-1)}f\|_{L^p(\Om)}
 \le e^{-t}\|f\|_{L^p(\Om)}
 \qquad(1\le p\le\infty,\ t\ge0),
\end{equation*}
and, whenever $1\le q\le p\le\infty$, there exists
$C=C(p,q,\Om)>0$ such that
\begin{equation*}
 \|e^{t(\Delta_N-1)}f\|_{L^p(\Om)}
 \le C e^{-t}
 \left(1+t^{-\frac32(\frac1q-\frac1p)}\right)
 \|f\|_{L^q(\Om)}
 \qquad(t>0).
\end{equation*}
\end{lemma}

We also use exponentially weighted maximal regularity for the damped
Neumann heat equation.

\begin{lemma}\label{lem:MR}
Let $r\in(1,\infty)$, and suppose that
\[
 z_t=\Delta z-z+f,\qquad \partial_\nu z=0,
 \qquad z(\cdot,0)=z_0\in W^{2,\infty}(\Om),
 \qquad\partial_\nu z_0=0.
\]
Then
\begin{equation*} \int_0^t e^{s-t}\|z(\cdot,s)\|_{W^{2,2}}^r\dd s
 \le C+C\int_0^t e^{s-t}\|f(\cdot,s)\|_{L^2}^r\dd s
\end{equation*}
for all times for which the solution exists.  The constant $C$ may depend
on $r,\Om$, and $z_0$, but is independent of $t$ and the terminal time of
the solution.
\end{lemma}

\begin{proof}
The function $e^{s/r}z(s)$ satisfies
\[
 \frac{\dd}{\dd s}\bigl(e^{s/r}z\bigr)
 -\Delta\bigl(e^{s/r}z\bigr)
 +\left(1-\frac1r\right)e^{s/r}z=e^{s/r}f.
\]
Because $r>1$, the Neumann realization
$-\Delta_N+1-1/r$ has a strictly positive zeroth-order part and admits
maximal $L^r$ regularity on $L^2(\Om)$.  Hence, on every
interval $(0,t)$, the assumption $z_0\in W^{2,\infty}(\Om)$ gives
\begin{align*}
 \int_0^t e^s\|z(s)\|_{W^{2,2}}^r\dd s
 &\le C+C\int_0^t e^s\|f(s)\|_{L^2}^r\dd s.
\end{align*}
The maximal-regularity constant is independent of the interval length.
Multiplication by $e^{-t}$ gives the asserted estimate; see also
\cite{CaoTao2021}.
\end{proof}

\section{The parabolic--elliptic--parabolic system for \texorpdfstring{$n=3$}{n=3}: Global boundedness}\label{sec:elliptic}


This section is devoted to the proof of Theorem~\ref{thm:main-elliptic-global}. We consider \eqref{sys} with $n=3$ and $\tau=0$ on $\Om\subset\R^3$, namely, the following parabolic--elliptic--parabolic system
\begin{equation}\label{sys-tau=0}
\left\{
\begin{aligned}
 u_t&=\nabla\!\cdot\!\bigl((1+u)^{m-1}\nabla u\bigr)
      -\chi\nabla\!\cdot(u\nabla v)-u-uw+\varphi,
      &&x\in\Om,\ t\in(0,\Tmax),\\
 0&=\Delta v-v+uw,
      &&x\in\Om,\ t\in(0,\Tmax),\\
 w_t&=\Delta w-w+v,
      &&x\in\Om,\ t\in(0,\Tmax),\\
 \partial_\nu u&=\partial_\nu v=\partial_\nu w=0,
      &&x\in\partial\Om,\ t\in(0,\Tmax),\\
 u(x,0)&=u_0(x),\
 w(x,0)=w_0(x),
      &&x\in\Om,
\end{aligned}
\right.
\end{equation}
and prove global existence for $m>4/3$.

\medskip

We first derive a weighted integral estimate for \(w\).

\begin{lemma}\label{lem:w-moments}
For every $b>1$, there exists $C_b>0$ such that the solution of \eqref{sys-tau=0}
satisfies
\begin{equation*} \int_0^t e^{s-t}\|w(\cdot,s)\|_{L^\infty}^b\dd s\le C_b
 \qquad t\in(0,\Tmax).
\end{equation*}
\end{lemma}

\begin{proof}
First, we claim that {$\|w(\cdot,t)\|_{L^2}$ is uniformly bounded in $(0,\Tmax)$.}
By Lemma \ref{lem:elliptic}  with $q=2$ and the
 equation satisfied by $v$ in (\ref{sys-tau=0}), we have
\[
 \|v(\cdot,s)\|_{L^2}
 \le C\|u(\cdot,s)w(\cdot,s)\|_{L^1}.
\]
This, together with Lemma \ref{lem:mass} and the spatial integrals of the
second and third equations in \eqref{sys-tau=0}, implies
$$
\int_0^t e^{s-t}\|v\|_{L^2}\dd s \le C\int_0^t e^{s-t}\int_\Om u w\dd x\dd s =C\int_\Om w(x,t)\dd x-Ce^{-t}\int_\Om w_0(x)\dd x \le C.
$$
Therefore,  the variation-of-constants formula and
Lemma \ref{lem:neumann-semigroup} yield
\begin{align*}
 \|w(\cdot,t)\|_{L^2}
 &\le e^{-t}\|w_0\|_{L^2}
   +\int_0^t e^{s-t}\|v(\cdot,s)\|_{L^2}\dd s\le C.
\end{align*}
The claim is proved.

It suffices to prove the asserted estimate for $b>4/3$.
Thanks to this claim and the second-order Gagliardo--Nirenberg inequality, we have
\begin{align}
 \|w\|_{L^\infty}
 &\le C\|w\|_{W^{2,2}}^{\frac34}
       \|w\|_{L^2}^{\frac14}+C\|w\|_{L^2}\le C\|w\|_{W^{2,2}}^{\frac34}+C.
 \label{w-GN}
\end{align}
Applying \eqref{w-GN} and Lemma~\ref{lem:MR} to the
 equation satisfied by $w$ in (\ref{sys-tau=0}), we obtain
$$
 \int_0^t e^{s-t}\|w\|_{L^\infty}^b\dd s
 \le C\int_0^t e^{s-t}\|w \|_{W^{2,2}}^{3b/4}\dd s+C
 \le C\int_0^t e^{s-t}\|v\|_{L^2}^{3b/4}\dd s+C.
$$
Moreover, by Lemmas \ref{lem:mass} and \ref{lem:elliptic},   and the
 equation satisfied by $v$ in (\ref{sys-tau=0}), it follows that
\begin{align*}
  \int_0^t e^{s-t}\|w\|_{L^\infty}^b\dd s
 &\le  C\int_0^t e^{s-t}\|uw\|_{L^1}^{3b/4}\dd s+C \\
 & \le  C\int_0^t e^{s-t}
 \|u\|_{L^1}^{3b/4}\|w\|_{L^\infty}^{3b/4}\dd s+C
 \le  C\int_0^t e^{s-t}\|w\|_{L^\infty}^{3b/4}\dd s+C.
\end{align*}
Since $3b/4<b$, Young's inequality absorbs the last integral on the right side and yields the desired estimate.
\end{proof}

We next derive uniform $L^p$ bounds for $u$.

\begin{proposition}\label{prop:u-Lp}
Let $(u,v,w)$ denote the solution to the system  \eqref{sys-tau=0}.  If $m>4/3$, then, for
every $p>1$,
\begin{equation*} \sup_{t\in(0,\Tmax)}\|u(\cdot,t)\|_{L^p}<\infty.
\end{equation*}
\end{proposition}

\begin{proof}
It suffices to take $p\ge2$. Testing the first equation in \eqref{sys-tau=0} with $pu^{p-1}$
and integrating by parts, we obtain
\begin{align}
 \frac{\dd}{\dd t}\int_\Om u^p
 &+p(p-1)\int_\Om
 u^{p-2}(1+u)^{m-1}|\nabla u|^2\notag\\
 &=-\chi(p-1)\int_\Om u^p\Delta v
   +p\int_\Om\varphi u^{p-1} -p\int_\Om u^p-p\int_\Om u^pw.
 \label{elliptic-u-exact}
\end{align}

By the Sobolev embedding theorem, on the bounded domain $\Om\subset\R^3$, we have
$$
W^{1,2}(\Omega)\hookrightarrow L^6(\Omega).
$$
Thus, by \(L^p\) interpolation, Young's inequality and Lemma~\ref{lem:mass}, we derive
\begin{align}
 \left\|u^{\frac{p+m-1}{2}}\right\|_{L^6}^2  &\leq C\left\|\nabla u^{\frac{p+m-1}{2}}\right\|_{L^2}^2
 +C\left(\int_\Om u\right)^{p+m-1} \notag\\
 &\le C\left\|\nabla u^{\frac{p+m-1}{2}}\right\|_{L^2}^2+C \notag\\
 &\le C\int_\Om
 u^{p-2}(1+u)^{m-1}|\nabla u|^2+C.
 \label{elliptic-coercivity}
\end{align}

Combining \eqref{elliptic-u-exact} and \eqref{elliptic-coercivity} and using
$w\ge0$, Young's inequality, and H\"older's inequality, we obtain
\begin{align*}
 &\frac{\dd}{\dd t}\int_\Om u^p+\int_\Om u^p
 +c\|u\|_{L^{3(p+m-1)}}^{p+m-1} \\
 &\le \frac{\dd}{\dd t}\int_\Om u^p+\int_\Om u^p
 +p(p-1)\int_\Om
 u^{p-2}(1+u)^{m-1}|\nabla u|^2+C \\
 &=-\chi(p-1)\int_\Om u^p\Delta v
   +p\int_\Om\varphi u^{p-1}
   -(p-1)\int_\Om u^p
   -p\int_\Om u^pw+C \\
 &\le \chi(p-1)\int_\Om u^{p+1} w
   +|\Om|\varphi_*^p+C,
\end{align*}
where the last inequality uses
$-\Delta v\le uw$ from the $v$ equation in \eqref{sys-tau=0}.
It follows that
\begin{align}
 \int_\Om u^p(x,t)
 +c\int_0^t e^{s-t}
 \|u(\cdot,s)\|_{L^{3(p+m-1)}}^{p+m-1}\dd s
 \le C\int_0^t e^{s-t}
 \|u\|_{L^{p+1}}^{p+1}\|w\|_{L^\infty}\dd s + C.
 \label{weighted-u}
\end{align}

By \(L^p\) interpolation and
Lemma~\ref{lem:mass},
\begin{equation}\label{u-interpolation-elliptic}
 \|u\|_{L^{p+1}}
 \le \|u\|_{L^{3(p+m-1)}}^{\gamma}
      \|u\|_{L^1}^{1-\gamma}
 \le C\|u\|_{L^{3(p+m-1)}}^{\gamma},
\end{equation}
where $\gamma\in(0,1)$ satisfies
\[
 \frac1{p+1}=\frac{\gamma}{3(p+m-1)}+1-\gamma.
\]
The existence of $\gamma\in(0,1)$ is  obvious since $1<p+1< 3(p+m-1)$.

Substituting \eqref{u-interpolation-elliptic} into \eqref{weighted-u}
and applying Young's inequality with any $b>1$, we obtain
\begin{align}\label{pf-prop-u-Lp}
 & \int_\Om u^p(x,t)
 +c\int_0^t e^{s-t}
 \|u(\cdot,s)\|_{L^{3(p+m-1)}}^{p+m-1}\dd s\cr
 \le\ & C\int_0^t e^{s-t}
 \|u\|_{L^{3(p+m-1)}}^{(p+1)\gamma}
 \|w\|_{L^\infty}\dd s+C\cr
 \le\ & C\int_0^t e^{s-t}
 \|u\|_{L^{3(p+m-1)}}^{\frac{(p+1)\gamma b}{b-1}}\dd s
 +C\int_0^t e^{s-t}\|w\|_{L^\infty}^b\dd s+C\cr
 \le\ & C\int_0^t e^{s-t}
 \left(\|u\|_{L^{3(p+m-1)}}^{p+m-1}\right)^{
 \frac{(p+1)\gamma b}{(p+m-1)(b-1)}}\dd s+C,
\end{align}
where the last inequality follows from Lemma \ref{lem:w-moments}.
The choice of   $\gamma$ shows that for $m>4/3$,
\begin{equation*}
 \frac{(p+1)\gamma b}{(p+m-1)(b-1)}
 ={b\over b-1}\frac{3p}{3p+3m-4}
 \longrightarrow\frac{3p}{3p+3m-4}<1
 \qquad\text{as }b\to\infty.
\end{equation*}
Thus we choose $b>1$  sufficiently large such that
$$
\frac{(p+1)\gamma b}{(p+m-1)(b-1)}<1,
$$
 and then it follows from (\ref{pf-prop-u-Lp}) and  Young's inequality that
\begin{align*}
 \int_\Om u^p(x,t)\dd x
 +c\int_0^t e^{s-t}
 \|u(\cdot,s)\|_{L^{3(p+m-1)}}^{p+m-1}\dd s\le \frac{c}{2}\int_0^t e^{s-t}
 \|u\|_{L^{3(p+m-1)}}^{p+m-1}\dd s+C.
\end{align*}
The desired conclusion follows immediately.
\end{proof}

The preceding estimates yield global boundedness for the
parabolic--elliptic--parabolic system.

\begin{proof}[Proof of Theorem \ref{thm:main-elliptic-global}]
The  $L^p$-estimates of $u$ in Proposition \ref{prop:u-Lp}, combined with the equations for $v$ and $w$ in the system (\ref{sys-tau=0}),
yield the integrability required for $u\nabla v$ in the Tao--Winkler
iteration lemma. The $L^\infty$-bootstrap result of Tao--Winkler
\cite[Lemma~A.1]{TaoWinkler2012JDE} and standard elliptic and
parabolic regularity theory yield the asserted boundedness and,
by Lemma~\ref{global_bounded_lemma}, global existence; see
\cite[Lemma~3.3]{wang2025} for this standard argument.
\end{proof}

\section{The parabolic--elliptic--parabolic system for \texorpdfstring{$n=3$}{n=3}: Finite-time blow-up}
\label{sec:blowup-subcritical}

Throughout this section, let $n=3$. We consider \eqref{sys-tau=0} on the ball $B_R(0)\subset\R^3$ with constant initial data $w_0\equiv W_0>0$, namely,
\begin{equation}\label{sys-radial-blowup}
\left\{
\begin{aligned}
 u_t&=\nabla\!\cdot\!\bigl((1+u)^{m-1}\nabla u\bigr)
      -\chi\nabla\!\cdot(u\nabla v)-u-uw+\varphi,
      &&x\in B_R(0),\ t\in(0,\Tmax),\\
 0&=\Delta v-v+uw,
      &&x\in B_R(0),\ t\in(0,\Tmax),\\
 w_t&=\Delta w-w+v,
      &&x\in B_R(0),\ t\in(0,\Tmax),\\
 \partial_\nu u&=\partial_\nu v=\partial_\nu w=0,
      &&x\in\partial B_R(0),\ t\in(0,\Tmax),\\
 u(x,0)&=u_0(x),\
 w(x,0)=W_0,
      &&x\in B_R(0).
\end{aligned}
\right.
\end{equation}
The following theorem immediately implies Theorem~\ref{thm:main-elliptic-blowup}.

\begin{theorem}
\label{thm:blowup-subcritical}
Consider \eqref{sys-radial-blowup} on $\Om=B_R(0)\subset\R^3$, and suppose
that $\varphi$ satisfies \eqref{source-ass} and is radial in space.
For every $\mu>0$, there exists $W^*=W^*(\mu)>0$ such that, if one of the following conditions holds:
\begin{enumerate}[label=\textup{(\roman*)}]
\item $1<m<4/3$ and $W_0>0$;
\item $m=4/3$ and $W_0>W^*$,
\end{enumerate}
then there exists a nonnegative radial $u_0\in W^{2,\infty}(\Om)$ with
$\partial_\nu u_0=0$ on $\partial\Om$ and $\int_\Om u_0=\mu$ such that
the corresponding solution of \eqref{sys-radial-blowup} blows up in finite time;
more precisely,
\[
 \Tmax<\infty,
 \qquad
 \limsup_{t\nearrow\Tmax}\|u(\cdot,t)\|_{L^\infty(\Om)}=\infty.
\]
\end{theorem}

To prove this theorem, we introduce the following J\"ager--Luckhaus transformation:
\begin{equation}\label{z-def}
 z(s,t):=3\int_0^{s^{1/3}}\rho^2u(\rho,t)\dd\rho,
 \qquad 0\le s\le R^3.
\end{equation}
A direct calculation yields
\begin{equation}\label{z-derivatives}
 z_s(s,t)=u(s^{1/3},t),\qquad
 z_{ss}(s,t)=\frac{u_r(s^{1/3},t)}{3s^{2/3}},
\end{equation}
and
\begin{equation}\label{z-boundary}
 z(0,t)=0,
 \qquad z(R^3,t)=\frac1{\Bunit}\int_\Om u(x,t)\dd x.
\end{equation}
Radial integration of the elliptic equation gives
\[
 r^2v_r(r,t)=\int_0^r\rho^2v(\rho,t)\dd\rho
             -\int_0^r\rho^2u(\rho,t)w(\rho,t)\dd\rho.
\]
Consequently, radial integration of the $u$ equation and the nonnegativity of $\varphi$ give
\begin{align}
 z_t\ge{}&9s^{4/3}(1+z_s)^{m-1}z_{ss}
 -3\chi z_s\int_0^{s^{1/3}}\rho^2v\dd\rho
 +3\chi z_s\int_0^{s^{1/3}}\rho^2uw\dd\rho\notag\\
 &-z-3\int_0^{s^{1/3}}\rho^2uw\dd\rho.
 \label{z-exact}
\end{align}

The proof of Theorem 4.1 is based on the following two lemmas. Lemma 4.2 establishes short-time bounds 
uniformly over initial data of prescribed mass, and uses these bounds to derive a scalar differential inequality for \(z\). Lemma 4.3 constructs a shrinking subsolution to this inequality. We then choose initial data of the prescribed mass to ensure the initial ordering, while the mass lower bound from Lemma 4.2 provides the required boundary ordering. By choosing the subsolution to concentrate at the origin within the uniform time interval, we can apply a comparison argument to deduce finite-time blow-up of \(u\).

\begin{lemma}\label{lem:short}
Fix $p\in(3/2,3)$, $\mu>0$ and $W_0>0$.
Then there exist $T_*>0$ and $K>0$, independent of the choice of $u_0$
with $\int_\Om u_0=\mu$, such that the solution $(u,v,w)$ of the problem (\ref{sys-radial-blowup})
satisfies
\begin{equation*} \frac\mu2\le\int_\Om u(\cdot,t)\le2\mu,
 \qquad \frac{W_0}{2}\le w(\cdot,t)\le2W_0,
 \qquad \|v(\cdot,t)\|_{L^p}\le K
\end{equation*}
for $0<t<\min\{T_*,\Tmax\}$.

In particular, given $\varepsilon\in(0,1/3)$, choose
$p=3/(1+3\varepsilon)\in(3/2,3)$ and set
\begin{equation*} B_0:=\frac{W_0\chi}{2}>0,
 \qquad C_0:=2W_0+1.
\end{equation*}
Then there exists $A_0=A_0(K)>0$ such that
\begin{equation*} z_t\ge \F[z]:=9s^{4/3}(1+z_s)^{m-1}z_{ss}
 -A_0s^{2/3-\varepsilon}z_s+B_0zz_s-C_0z
\end{equation*}
on $(0,R^3)\times(0,\min\{T_*,\Tmax\})$.
\end{lemma}

\begin{proof}
The mass identity gives
\[
 \frac{\dd}{\dd t}\int_\Om u
 \le-\int_\Om u+|\Om|\varphi_*,
\]
hence
\begin{align*}
 \int_\Om u(\cdot,t)
 \le\mu e^{-t}+|\Om|\varphi_*(1-e^{-t})\le\mu+|\Om|\varphi_*t\le\frac{3\mu}{2}
\end{align*}
for $0<t<\min\{\mu/(2|\Om|\varphi_*),\Tmax\}$.

Lemma~\ref{lem:elliptic} gives
\[
 \|v\|_{L^p}\le C_p\|uw\|_{L^1}.
\]
The variation-of-constants formula for $w$ reads
\[
 w(\cdot,t)=W_0e^{-t}
 +\int_0^t e^{(t-s)(\Delta_N-1)}v(\cdot,s)\dd s.
\]
Using Lemma~\ref{lem:neumann-semigroup}, the preceding estimate for $v$,
and the mass bound above, we obtain
\begin{align*}
 \|w(\cdot,t)\|_{L^\infty}
 &\le W_0e^{-t}+C\int_0^t e^{s-t}
 \bigl(1+(t-s)^{-3/(2p)}\bigr)
 \|v(\cdot,s)\|_{L^p}\dd s\\
 &\le W_0+C\mu\int_0^t
 \bigl(1+(t-s)^{-3/(2p)}\bigr)
 \|w(\cdot,s)\|_{L^\infty}\dd s.
\end{align*}
Set $M(t):=\sup_{0\le s\le t}\|w(\cdot,s)\|_{L^\infty}$.
Since $3/(2p)<1$, we may choose $T_*>0$, with
$T_*\le\mu/(2|\Om|\varphi_*)$, so small that
\[
 C\mu\left(T_*+\frac{T_*^{1-3/(2p)}}{1-3/(2p)}\right)\le\frac12.
\]
For $0<t<\min\{T_*,\Tmax\}$, taking the supremum in the preceding
estimate gives
\[
 M(t)\le W_0+C\mu\left(t+\frac{t^{1-3/(2p)}}{1-3/(2p)}\right)M(t)
 \le W_0+\frac12M(t),
\]
and hence $M(t)\le2W_0$.

Using this upper bound in the mass identity, we obtain
\[
 \frac{\dd}{\dd t}\int_\Om u
 \ge-(1+2W_0)\int_\Om u
\]
and therefore
\[
 \int_\Om u(\cdot,t)\ge\mu e^{-(1+2W_0)t}\ge\frac\mu2
\]
after decreasing $T_*$ so that $T_*\le\log 2/(1+2W_0)$.
Moreover, Lemma~\ref{lem:elliptic} and the upper bounds for the mass of
$u$ and for $w$ yield
\[
 \|v(\cdot,t)\|_{L^p}\le C_p\|uw\|_{L^1}
 \le3C_p\mu W_0=:K.
\]

Finally, the nonnegativity of $v$ and parabolic comparison give
\[
 w(\cdot,t)\ge W_0e^{-t}\ge\frac{W_0}{2}
\]
for the same time interval, since $T_*\le\log 2$.
Both $T_*$ and $K$ are independent of the choice of $u_0$ with the prescribed mass.

Now let $\varepsilon\in(0,1/3)$ and take $p=3/(1+3\varepsilon)$.
For $t\in(0,\min\{T_*,\Tmax\})$, the preceding bounds and
H\"older's inequality show that
\[
 \int_0^{s^{1/3}}\rho^2v(\rho,t)\dd\rho
 \le CKs^{1-1/p}=CKs^{2/3-\varepsilon}.
\]
Moreover, the bounds for $w$ give
\[
 3\chi z_s\int_0^{s^{1/3}}\rho^2uw\dd\rho
 \ge\frac{W_0\chi}{2}zz_s,
 \qquad
 -3\int_0^{s^{1/3}}\rho^2uw\dd\rho\ge-2W_0z.
\]
Substituting these estimates into \eqref{z-exact} yields the asserted scalar inequality.
\end{proof}

Fix $\varepsilon\in(0,1/3)$ and let $\F$ be the operator defined in Lemma~\ref{lem:short}.
Its coefficients $A_0,B_0,C_0$ depend on $W_0$ but are independent of the
choice of $u_0$ with prescribed mass. We seek a shrinking subsolution of the form
\begin{equation}\label{barrier}
 \uz(s,t)=a(t)f\left(\frac{s}{b(t)}\right),
 \qquad 0\le t<T_b,
\end{equation}
where
\[
\begin{aligned}
 a(t)&=a_0e^{-\Lambda t},\\
 b(t)&=\bigl(b_0^{1-\sigma}-\kappa_0(1-\sigma)t\bigr)^{1/(1-\sigma)},
 \qquad T_b:=\frac{b_0^{1-\sigma}}{\kappa_0(1-\sigma)},
\end{aligned}
\]
and
\begin{equation}\label{profile}
f(\xi)=\begin{cases}(k+1)(k+2)\displaystyle\int_0^\xi \eta(1-\eta)^k\dd \eta, ~&\xi\in[0,1],\\
1, ~&\xi>1
\end{cases},
\qquad k=\frac{m}{m-1}.
\end{equation}
The parameters $a_0,b_0,\kappa_0,\Lambda>0$ and $\sigma\in(1/3,1)$
will be chosen below.

\begin{lemma}\label{lem:subcritical-barrier}
For every $a_0>0$, there exists $B^\ast=B^\ast(a_0)>0$ such that, if one of the following conditions holds:
\begin{enumerate}[label=\textup{(\roman*)}]
\item $1<m<4/3$ and $B_0>0$;
\item $m=4/3$ and $B_0>B^\ast$,
\end{enumerate}
then, for any $\sigma\in(1/3,1)$, $\kappa_0>0$ and $\Lambda>C_0$,
the function $\uz$ defined in \eqref{barrier} satisfies
\[
\uz_t\le\F[\uz]
\qquad\text{in }(0,\infty)\times(0,T_b)
\]
for every sufficiently small $b_0>0$.
\end{lemma}

\begin{proof}
Set $\xi=s/b(t)$. Using $a'=-\Lambda a$ and $b'=-\kappa_0b^\sigma$,
a direct calculation yields
\begin{align}
\uz_t-\F[\uz]={}&-9ab^{-2/3}\xi^{4/3}
\left(1+\frac abf'\right)^{m-1}f''
+A_0ab^{-1/3-\varepsilon}\xi^{2/3-\varepsilon}f'
\notag\\
&-B_0a^2b^{-1}ff'-(\Lambda-C_0)af
+\kappa_0ab^{\sigma-1}\xi f'.
\label{residual}
\end{align}

\noindent\textbullet\ The case $\xi\in(1,\infty)$.

Since $\Lambda>C_0$, \eqref{residual} gives
\[
\uz_t-\F[\uz]=-(\Lambda-C_0)af<0.
\]

\noindent\textbullet\ The case
$\xi\in[0,{1}/{(2k+2)}]$.

For such $\xi$,
\[
f'(\xi)\le (k+1)(k+2)\xi,
\]
and
\[
f''(\xi)
\ge
(k+1)(k+2)
\left(1-\frac12\right)
\left(1-\frac1{2(k+1)}\right)^{k-1}
\ge
\frac{(k+1)(k+2)}{2^k}.
\]

Substituting these estimates into \eqref{residual} and using
$\Lambda>C_0$, we obtain
\begin{align*}
\uz_t-\F[\uz]
&\le
(k+1)(k+2)a
\Big[
-9b^{-2/3}\xi^{4/3}2^{-k}
+A_0 b^{-1/3-\varepsilon}\xi^{5/3-\varepsilon}
+\kappa_0b^{\sigma-1}\xi^2
\Big]
\\
&=
(k+1)(k+2)ab^{-2/3}\xi^{4/3}
\Big[
-\frac9{2^k}
+A_0 b^{1/3-\varepsilon}\xi^{1/3-\varepsilon}
+\kappa_0b^{\sigma-\frac13}\xi^{\frac23}
\Big]
\\
&\le
(k+1)(k+2)ab^{-2/3}\xi^{4/3}
\Big[
-\frac9{2^k}
+A_0 b_0^{1/3-\varepsilon}
+\kappa_0b_0^{\sigma-\frac13}
\Big]
\\
&\le0,
\end{align*}
provided that
\begin{equation}\label{cond1}
\boxed{
A_0 b_0^{1/3-\varepsilon}
+\kappa_0b_0^{\sigma-\frac13}
\le\frac9{2^k}
}.
\end{equation}

\noindent\textbullet\ The case
$\xi\in[{1}/{(2k+2)},1]$.

In this case,
\[
f(\xi)
\ge
(k+1)(k+2)
\int_0^{\frac1{2k+2}}
\eta(1-\eta)^k\,d\eta
=:(k+1)(k+2)c_k^\ast>0,
\]
and
\[
\frac{k+2}{2}(1-\xi)^k
\le
f'(\xi)
\le
(k+1)(k+2)(1-\xi)^k,
\]
as well as
\[
|f''(\xi)|
\le
(k+1)^2(k+2)(1-\xi)^{k-1}.
\]

Substituting these estimates into \eqref{residual}, we obtain
\begin{align}
\uz_t-\F[\uz]
\le{}&
(k+1)(k+2)ab^{-1}
\Big[
9b^{1/3}
\left(1+\frac abf'\right)^{m-1}
(k+1)(1-\xi)^{k-1}
\notag\\
&\quad
-\frac{B_0ac_k^\ast(k+2)}2(1-\xi)^k
-(\Lambda-C_0)bc_k^\ast
\notag\\
&\quad
+A_0 b^{2/3-\varepsilon}(1-\xi)^k
+\kappa_0b^\sigma(1-\xi)^k
\Big],
\label{residua2}
\end{align}
where we used $\xi\le1$.

Choose $b_0$ sufficiently small so that
\begin{equation}\label{Tb}
\boxed{
\frac{\Lambda b_0^{1-\sigma}}
{\kappa_0(1-\sigma)}
=
\Lambda T_b
\le\log2
}
\end{equation}
and
\begin{equation}
\label{smallb1}
\boxed{
A_0 b_0^{2/3-\varepsilon}
+\kappa_0b_0^\sigma
\le
\frac{B_0a_0c_k^\ast(k+2)}8
}.
\end{equation}

Since \eqref{Tb} implies
\[
\frac{a_0}{2}\le a(t)\le a_0,
\]
it follows from \eqref{smallb1} that
\[
A_0 b^{2/3-\varepsilon}(1-\xi)^k
+\kappa_0b^\sigma(1-\xi)^k
\le
\frac{B_0ac_k^\ast(k+2)}4(1-\xi)^k.
\]

Using this in \eqref{residua2}, we obtain
\begin{align}
\uz_t-\F[\uz]
\le{}&
(k+1)(k+2)ab^{-1}
\Big[
9b^{1/3}
\left(
1+\frac{(k+1)(k+2)a_0}{b}(1-\xi)^k
\right)^{m-1}
\notag\\
&\qquad\qquad\qquad\qquad\qquad\qquad
\times(k+1)(1-\xi)^{k-1}
\notag\\
&\quad
-\frac{B_0ac_k^\ast(k+2)}4(1-\xi)^k
-(\Lambda-C_0)c_k^\ast b
\Big].
\label{residua3}
\end{align}

To show that the right-hand side is nonpositive, we distinguish the
cases $(1-\xi)^k\le b$ and $(1-\xi)^k>b$.

If $(1-\xi)^k\le b$, then
\[
b^{1/3}(1-\xi)^{k-1}
\le
b^{\frac13+\frac{k-1}{k}}
=
b^{\frac43-\frac1k}.
\]
Let
\[
M_0:=1+(k+1)(k+2)a_0.
\]
Then
\begin{align}
&9(k+1)b^{1/3}
\left(
1+\frac{(k+1)(k+2)a_0}{b}(1-\xi)^k
\right)^{m-1}
(1-\xi)^{k-1}
\notag\\
&\qquad
\le
9(k+1)b^{\frac43-\frac1k}M_0^{m-1}
\notag\\
&\qquad
\le
9(k+1)b_0^{\frac13-\frac1k}
M_0^{m-1}b
\notag\\
&\qquad
\le
(\Lambda-C_0)c_k^\ast b,
\label{f''est1}
\end{align}
provided that
\begin{equation}
\label{smallb-case1}
\boxed{
9(k+1)M_0^{m-1}
b_0^{\frac13-\frac1k}
\le
(\Lambda-C_0)c_k^\ast
}.
\end{equation}

Thus \eqref{f''est1} implies that the right-hand side of
\eqref{residua3} is nonpositive.

Now suppose that $(1-\xi)^k>b$. Since
\[
\frac{(1-\xi)^k}{b}>1,
\]
we have
\[
1+\frac{(k+1)(k+2)a_0}{b}(1-\xi)^k
\le
\frac{M_0(1-\xi)^k}{b}.
\]
Hence
\begin{align}
&9(k+1)b^{1/3}
\left(
1+\frac{(k+1)(k+2)a_0}{b}(1-\xi)^k
\right)^{m-1}
(1-\xi)^{k-1}
\notag\\
&\qquad
\le
9(k+1)M_0^{m-1}
b^{\frac43-m}
(1-\xi)^{km-1}
\notag\\
&\qquad
\le
9(k+1)M_0^{m-1}
b_0^{\frac43-m}
(1-\xi)^k
\notag\\
&\qquad
\le
\frac{B_0a_0c_k^\ast(k+2)}8
(1-\xi)^k
\notag\\
&\qquad
\le
\frac{B_0ac_k^\ast(k+2)}4
(1-\xi)^k,
\label{f''est2}
\end{align}
provided that
\begin{equation}
\label{smallb2}
\boxed{
9(k+1)M_0^{m-1}
b_0^{\frac43-m}
\le
\frac{B_0a_0c_k^\ast(k+2)}8
}.
\end{equation}

Here we used $km-1-k=m-1>0$.

Thus \eqref{f''est2} implies that the right-hand side of
\eqref{residua3} is nonpositive.

Finally, we verify that all the boxed conditions can be satisfied
simultaneously. Since
\[
\frac13-\varepsilon>0,\qquad
\sigma-\frac13>0,\qquad
1-\sigma>0,\qquad
\frac13-\frac1k>0,
\]
the conditions \eqref{cond1}, \eqref{Tb}, \eqref{smallb1}, and
\eqref{smallb-case1} hold for all sufficiently small $b_0>0$.

If $1<m<\frac43$, then
\[
\frac43-m>0,
\]
so \eqref{smallb2} also holds after decreasing $b_0$ if necessary.

If $m=\frac43$, then $k=4$ and \eqref{smallb2} reduces to
\[
9(k+1)M_0^{1/3}
\le
\frac{B_0a_0c_k^\ast(k+2)}8.
\]
This condition holds whenever $B_0>B^\ast(a_0)$ for a suitable constant
$B^\ast(a_0)>0$.

Thus all the required conditions can be satisfied simultaneously by a
suitable choice of $b_0>0$. This completes the proof.
\end{proof}

We now prove Theorem~\ref{thm:blowup-subcritical} using the preceding results.

\begin{proof}[Proof of Theorem~\ref{thm:blowup-subcritical}]

Fix a prescribed mass $\mu>0$ and the parameters in Lemma~\ref{lem:short}, and set
$a_0=\mu/(2\Bunit)$. At $m=4/3$, choose
\[
 W^*(\mu):=\frac{2B^\ast(a_0)}{\chi},
\]
where $B^\ast$ is the threshold from Lemma~\ref{lem:subcritical-barrier}.
Since $B_0=\chi W_0/2$, the required condition $B_0>B^\ast(a_0)$ is
equivalent to $W_0>W^*(\mu)$.

Fix $\sigma\in(1/3,1)$,
$\kappa_0>0$, and $\Lambda>C_0$.  Lemma~\ref{lem:subcritical-barrier}
allows $b_0\in(0,R^3)$ to be chosen
so small that $T_b<T_*$ and $\uz_t\le\F[\uz]$ on
$(0,R^3)\times(0,T_b)$.  Define
\[
 u_0(x)=\frac{a_0}{b_0}
 f'\left(\frac{|x|^3}{b_0}\right)+\delta,
\]
where
\[
\delta
=
\frac{\mu}{\Bunit R^3}
-
\frac{a_0}{R^3}.
\]
Because $b_0<R^3$ and $f(R^3/b_0)=1$, this choice gives
$\delta=\mu/(2|\Om|)>0$ and $\int_\Om u_0=\mu$.
Moreover, $u_0\in W^{2,\infty}(\Om)$ is nonnegative, radial, and
constant near $\partial\Om$, so it satisfies the Neumann condition.

Indeed, $z(s,0)=a_0f(s/b_0)+\delta s\ge\uz(s,0)$, while
Lemma~\ref{lem:short} yields
\[
 z(R^3,t)\ge\frac{\mu}{2\Bunit}
 =a_0\ge a(t)=\uz(R^3,t)
 \qquad(0<t<\min\{T_b,\Tmax\}).
\]
Since $z(0,t)=\uz(0,t)=0$, the ordering $\uz\le z$ holds on the
parabolic boundary. Fix $0<T_0<\min\{T_b,\Tmax\}$, let $Z=\uz-z$, and
choose $\lambda>0$ such that $\lambda>B_0z_s-C_0$ on
$[0,R^3]\times[0,T_0]$; this is possible because $z_s=u$ is bounded there.
If $e^{-\lambda t}Z$ had a positive maximum, it would occur at some
$(s_0,t_0)\in(0,R^3)\times(0,T_0]$. At this point,
$Z_s=0$ and $Z_{ss}\le0$, so the differential inequalities for $z$ and
$\uz$ imply
\[
 \lambda Z\le Z_t
 \le9s^{4/3}(1+z_s)^{m-1}Z_{ss}+(B_0z_s-C_0)Z
 <\lambda Z,
\]
a contradiction. Thus $\uz\le z$ on $[0,R^3]\times[0,T_0]$.
Since $T_0$ is arbitrary, we obtain
\begin{equation}\label{final-comparison}
 z(s,t)\ge\uz(s,t)
 \quad(0\le s\le R^3,\ 0\le t<\min\{T_b,\Tmax\}).
\end{equation}
At $s=b(t)/2$,
\[
 z\left(\frac{b(t)}2,t\right)
 \le\frac{b(t)}2\|u(\cdot,t)\|_{L^\infty},
\]
and hence \eqref{final-comparison} implies
\begin{equation}\label{lower-u-infty}
 \|u(\cdot,t)\|_{L^\infty}
 \ge\frac{2a(t)f(1/2)}{b(t)}.
\end{equation}
Recalling $a(t)=a_0e^{-\Lambda t}$ and
$b(t)=[\kappa_0(1-\sigma)(T_b-t)]^{1/(1-\sigma)}$, we have
\[
 \frac{2a(t)f(1/2)}{b(t)}
 =\frac{2a_0e^{-\Lambda t}f(1/2)}
 {[\kappa_0(1-\sigma)(T_b-t)]^{1/(1-\sigma)}}
 \longrightarrow\infty\quad\text{as }t\nearrow T_b.
\]
Hence $\Tmax\le T_b<\infty$.  If $\Tmax=T_b$, this directly gives blow-up of
$u$.  If $\Tmax<T_b$, Lemma~\ref{lem:short} keeps $w$ bounded up to
$\Tmax$, and Lemma~\ref{global_bounded_lemma} forces $u$ to become unbounded.
\end{proof}

\begin{remark}
\label{rem:nonconstant-w0}
The proof of Theorem~\ref{thm:blowup-subcritical} also applies to radial $w_0\in W^{2,\infty}(\Om)$ satisfying
$\partial_\nu w_0=0$ and
\[
 0<W_-\le w_0(x)\le W_+
 \qquad(x\in\overline\Om).
\]
In this case, the bounds for $w$ in Lemma~\ref{lem:short} become
$W_-/2\le w(\cdot,t)\le2W_+$, and $T_*$ and $K$ may depend on
$W_-$ and $W_+$ but remain independent of the choice of $u_0$ with
$\int_\Om u_0=\mu$. The scalar inequality in Lemma~\ref{lem:short} then holds with
$B_0=\chi W_-/2$ and $C_0=2W_++1$.  Thus the blow-up construction
works for every $W_->0$ when $m<4/3$, and for sufficiently large $W_-$
when $m=4/3$; the condition $W_->W^*(\mu)$ in
Theorem~\ref{thm:blowup-subcritical} suffices.
Constancy of $w_0$ merely makes the coefficient condition
more transparent.
\end{remark}

\section{The parabolic--elliptic--parabolic system for \texorpdfstring{$n\ge4$}{n >= 4}: Finite-time blow-up}\label{sec:high-blowup}

The following proposition refines Theorem~\ref{thm:main-high-blowup}.
After establishing common estimates, we first consider
$0<m\le1$ using a linear moment and then $1<m\le2-2/n$ using a logarithmic moment.

\begin{proposition}
\label{prop:high-blowup-restatement}
Let $n\ge4$, $\Om=B_R(0)$, $\tau=0$, and $0<m\le2-2/n$.
Suppose that \eqref{source-ass} holds and that $\varphi$ is radial in space.
Fix $\mu>0$. If $m<2-2/n$, take any $W>0$; if $m=2-2/n$, assume
\begin{equation*} W>\frac{3e(\widetilde{C}_3+1)}{\widetilde{C}_1C_\mu^2},
 \qquad C_\mu=\frac12\log\left(1+\frac{\mu}{\Bunit}\right),
\end{equation*}
where $\widetilde{C}_1$ and $\widetilde{C}_3$ are the constants in Proposition~\ref{n ge 4}
at $m=2-2/n$, with $M=\max\{\mu/\Bunit,R^n\varphi_*\}$,
and are independent of $W$.
Then there exists a smooth nonnegative radial $u_0$ satisfying
\eqref{data-ass} and $\int_\Om u_0=\mu$ such that the solution of
\eqref{sys} with $w_0\equiv W$ satisfies
\[
 \Tmax<\infty,\qquad
 \limsup_{t\nearrow\Tmax}\|u(\cdot,t)\|_{L^\infty(\Om)}=\infty.
\]
\end{proposition}

For the moment, let $u_0$ be any smooth nonnegative radial function with
$\int_\Om u_0=\mu>0$, and take $w_0\equiv W>0$.
Recall that integrating the first equation gives
\[
 \frac{\dd}{\dd t}\int_\Om u\le-\int_\Om u+|\Om|\varphi_*.
\]
Hence
\begin{equation}\label{high-mass}
 \int_\Om u(\cdot,t)\le\Bunit M
 \qquad(0<t<\Tmax),
\end{equation}
where \[
 M=\max\left\{\frac{\mu}{\Bunit},R^n\varphi_*\right\}.
\]
In particular, $M$ is independent of $W$ and the concentration scale.

We now consider \eqref{sys} with $\tau=0$ in
$\Om=B_R(0)\subset\R^n$, where $n\ge4$ and $0<m\le2-2/n$.
The source $\varphi$ satisfies \eqref{source-ass} and is radial in space.
Unlike in the three-dimensional case, for $n\ge4$ it is difficult to obtain a uniform upper bound for $w$, which prevents a direct use of the reduction introduced by Tao and Winkler \cite{TaoWinkler2021}. This is also why their blow-up result for $m=1$ was restricted to $n=2,3$, and why the later extension by Wang and Wang \cite{wang2025} to $m<1$ remained confined to the same dimensions. To overcome this difficulty, we instead establish a weighted radial monotonicity property for $w$ (note that radial monotonicity of the initial data does not in general imply radial monotonicity of the solution), which provides the control needed to simplify the equation.

\begin{lemma}\label{lem:w-weight decreasing}
Let $w(r,0)\equiv W>0$. For $0<t<\Tmax$, $e^{-r}w(r,t)$ is nonincreasing in $r$. Thus, for $0\le r_1\le r_2\le R$,
\begin{equation*}
 w(r_2,t)\le e^{r_2-r_1} w(r_1,t),
\end{equation*}
and, in particular,
\begin{equation*} w(r_2,t)\le e^{R} w(r_1,t).
\end{equation*}
\end{lemma}

\begin{proof}
The variation-of-constants formula gives
$$w(\cdot,t)=We^{-t}+\int_0^t
 e^{(t-\sigma)(\Delta_N-1)}v(\sigma)\dd\sigma.$$

Combining this with the elliptic equation $(1-\Delta_N)v=uw$ gives
\begin{equation}\label{high-positive-resolvent}
 (1-\Delta_N)w(t)
 =We^{-t}+\int_0^t
 e^{(t-\sigma)(\Delta_N-1)}(uw)(\sigma)\dd\sigma\ge0.
\end{equation}

In radial coordinates, this reads
$$\left(r^{n-1}w_r\right)_r\le r^{n-1}w.$$
Fix $t\in(0,\Tmax)$. On any connected component of $\{r\in(0,R):w_r(r,t)>0\}$ with left endpoint $r_0$, multiplying by $2r^{n-1}w_r$ gives
\begin{align*}
\left(r^{2(n-1)}w^2_r\right)_r&\le r^{2(n-1)}\left(w^2\right)_r\\
&=\left(r^{2(n-1)}w^2\right)_r-2(n-1)r^{2n-3}w^2\\
&\le \left(r^{2(n-1)}w^2\right)_r.
\end{align*}

Integrating over $(r_0,r)$ and using $w_r(r_0,t)=0$, also when $r_0=0$, we obtain
\begin{align*}
r^{2(n-1)}w_r^2(r,t)
&\le r^{2(n-1)}w^2(r,t)-r_0^{2(n-1)}w^2(r_0,t)\\
&\le r^{2(n-1)}w^2(r,t).
\end{align*}
Thus $w_r\le w$ wherever $w_r>0$, and the same inequality is immediate where $w_r\le0$. Hence $e^{-r}w$ is nonincreasing in $r$.
\end{proof}

Let $$z(s,t)=n\int_0^{s^{1/n}}\rho^{n-1}u(\rho,t)\dd\rho$$ and $$Q(s,t)=n\int_0^{s^{1/n}}\rho^{n-1}uw(\rho,t)\dd\rho.$$

Then
\begin{equation}\label{fund est}
 z_s=u(s^{1/n},t)\ge0,\qquad
 Q_s=z_sw(s^{1/n},t),\qquad
 0\le z(s,t)\le M.
\end{equation}

The above lemma yields the following simple but useful consequence.
\begin{corollary}\label{cor:radial-interaction}
For $0<s\le R^n$, the functions $z$ and $Q$ defined above satisfy the following estimate.
\begin{equation*}\left[z(s,t)Q(s,t)\right]_s\le (1+e^R) z_s(s,t)Q(s,t).
\end{equation*}
\end{corollary}
\begin{proof}
Using Lemma~\ref{lem:w-weight decreasing} and repeatedly applying \eqref{fund est}, we obtain
\begin{align*}
 z(s,t)Q_s(s,t)&=\int_0^s z_s(\sigma,t)\dd \sigma\cdot z_s(s,t)w(s^{\frac1n},t)\\
&\le e^R \int_0^s z_s(\sigma,t)w(\sigma^{\frac1n},t)\dd \sigma\cdot z_s(s,t)\\
&\le e^R \int_0^s Q_s(\sigma,t)\dd \sigma\cdot z_s(s,t)\\
&\le e^{R}Q(s,t)z_s(s,t).
\end{align*}

Adding $Q(s,t)z_s(s,t)$ to both sides proves the corollary.
\end{proof}

Radial integration of \eqref{sys} gives
\begin{align}
 z_t={}&n^2s^{2-\frac{2}{n}}(1+z_s)^{m-1}z_{ss}
 +\chi z_s\left(Q-n\int_0^{s^{1/n}}\rho^{n-1}v\dd\rho\right)\notag\\
 &-z-Q+n\int_0^{s^{1/n}}\rho^{n-1}\varphi\dd\rho.
 \label{high-z-equation}
\end{align}

\subsection{Proof of Proposition~\ref{prop:high-blowup-restatement}: The case~$0<m\le1$}

For $0<m\le1$, the mass bound controls the diffusion term directly.

\begin{proof}[Proof of Proposition~\ref{prop:high-blowup-restatement} for $0<m\le1$]
Fix $\mu>0$ and $W>0$. For $0<s_0<\min\{1,R^n\}$ to be chosen below,
take a smooth nonnegative radial $u_0$ satisfying \eqref{data-ass},
supported in $B_{(s_0/2)^{1/n}}(0)$, with $\int_\Om u_0=\mu$,
and set $w_0\equiv W$.
Set
\begin{equation}\label{high-linear-moment}
 J_1(t)=\int_0^{s_0}z(s,t)\dd s
 =n\int_0^{s_0^{1/n}}r^{n-1}(s_0-r^n)u(r,t)\dd r.
\end{equation}
Suppose, to the contrary, that the solution is global. Then
\begin{equation}\label{J_1(0) low}
 J_1(0)\ge\frac{\mu}{2\Bunit}s_0,\qquad
 z(s_0,t)\ge\frac{J_1(t)}{s_0},~ t\ge 0.
\end{equation}

Set \[
 T_1^*:=\sup\left\{T\le1\ \middle|\
 J_1(t)\ge\frac{\mu}{3\Bunit}s_0,~~ 0\le t<T
 \right\}.
\]
We work on $[0,T_1^*)$.

Integrating \eqref{high-z-equation} over $(0,s_0)$ gives

\begin{align}\notag
 J'_1(t)&=\left(\int_0^{s_0}z(s,t)\dd s\right)_t\\
&\ge n^2\int_0^{s_0}s^{2-\frac2n}(1+z_s)^{m-1}z_{ss}\dd s+\chi\int_0^{s_0}z_s
\left(Q-n\int_0^{s^{1/n}}\rho^{n-1}v\dd\rho\right)\dd s\notag\\
&\quad-\int_0^{s_0}(z+Q)\dd s\notag\\
&=:I+II+III,
\label{high-linear-decomposition}
\end{align}
where we have dropped the nonnegative source term.

Since $(1+z_s)^m\le1+z_s$ for $0<m\le1$, integration by parts gives
\begin{align*}
 I&=\frac{n^2}{m}\int_0^{s_0}s^{2-\frac2n}\dd \left((1+z_s)^{m}\right)\\
 &=\frac{n^2}{m}s_0^{2-\frac2n}(1+z_s(s_0,t))^m
 -\frac{2n(n-1)}m\int_0^{s_0}s^{1-\frac2n}(1+z_s)^m\dd s\\
 &\ge-\frac{2n(n-1)}m s_0^{1-\frac2n}
                  \int_0^{s_0}(1+z_s)\dd s\\
 &\ge-\frac{2n(n-1)}m(1+M)s_0^{1-\frac2n}.
\end{align*}

 To estimate ${II}$, we write it as

\begin{align*}
{II}&=\chi\int_0^{s_0}z_sQ\dd s
-n\chi\int_0^{s_0}z_s
\left(\int_0^{s^{1/n}}\rho^{n-1}v\dd\rho\right)\dd s\\
&=:{II}_1-{II}_2.
\end{align*}

Direct computation gives
\begin{align}\notag
{II}_1&\ge \chi \int_0^{s_0}z_s(s,t)Q(s,t)\dd s\\\notag
&\ge \chi(1+e^R)^{-1}  \int_0^{s_0}[z(s,t)Q(s,t)]_s\dd s\\
&= \chi(1+e^R)^{-1}  z(s_0,t)Q(s_0,t)=:C_1 z(s_0,t)Q(s_0,t),\label{high-linear-II1-est}
\end{align}
where the second inequality follows from Corollary~\ref{cor:radial-interaction}.

We next estimate ${II}_2$. By the Brezis--Strauss estimate \cite{BrezisStrauss1973},
\[
 \|v\|_{W^{1,p}}\le C_p\|uw\|_{L^1(B_R(0))},~~ 1\le p<\frac{n}{n-1}.
\]
Taking $p=1$ and applying Sobolev embedding gives
\[
 \|v\|_{L^{n/(n-1)}(B_R(0))}\le C \|v\|_{W^{1,1}(B_R(0))}\le C\|uw\|_{L^1(B_R(0))}.
\]
Hence
\begin{align*}
 n\int_0^{s^{1/n}}\rho^{n-1}v(\rho,t)\dd\rho
 &\le n\left(\int_0^{s^{1/n}}\rho^{n-1}v^{\frac{n}{n-1}}(\rho,t)\dd\rho\right)^{\frac{n-1}{n}}\left(\int_0^{s^{1/n}}\rho^{n-1}\dd\rho\right)^{\frac{1}{n}}\\
&\le C s^{1/n}\|uw\|_{L^1(B_R(0))}\\
&\le Cs^{1/n} Q(R^n,t),~~0<s<R^n.
\end{align*}

By Lemma~\ref{lem:w-weight decreasing},
\begin{align*}
Q(R^n,t)-Q(s_0,t)&=n\int_{s_0^{\frac{1}{n}}}^{R}\rho^{n-1}uw\dd\rho\\
&\le ne^Rw(s_0^{\frac{1}{n}},t)\int_{s_0^{\frac{1}{n}}}^{R}\rho^{n-1}u\dd\rho\\
&\le e^R M w(s_0^{\frac{1}{n}},t).
\end{align*}
Since $t\in[0,T_1^*)$ implies $z(s_0,t)\ge J_1(t)/s_0\ge\mu/(3\Bunit)$,
Lemma~\ref{lem:w-weight decreasing} also gives
\begin{align*}
Q(s_0,t)&=n\int^{s_0^{\frac{1}{n}}}_0\rho^{n-1}uw(\rho,t)\dd\rho\\
&\ge ne^{-R}w(s_0^{\frac{1}{n}},t)\int^{s_0^{\frac{1}{n}}}_0\rho^{n-1}u\dd\rho\\
&\ge e^{-R}w(s_0^{\frac{1}{n}},t)z(s_0,t)\\
&\ge e^{-R}w(s_0^{\frac{1}{n}},t)\cdot\frac{\mu}{3\Bunit}.
\end{align*}

Combining these estimates gives
\begin{align*}
Q(R^n,t)\le \left(1+e^{2R}\frac{3\Bunit M}{\mu}\right)Q(s_0,t).
\end{align*}
Using this in the estimate for $v$ gives
\begin{align*}
 n\int_0^{s^{1/n}}\rho^{n-1}v(\rho,t)\dd\rho\le Cs^{1/n} Q(s_0,t),~~0<s<R^n,
\end{align*}
and hence
\begin{align}\notag
{II}_2:&= \chi\int_0^{s_0}z_s
 \left(n\int_0^{s^{1/n}}\rho^{n-1}v\dd\rho\right)\dd s\\\notag
&\le C s_0^{\frac{1}{n}} Q(s_0,t)\int_0^{s_0}z_s\dd s\\\notag
&\le C M s_0^{1/n}Q(s_0,t)=:C_2 s_0^{1/n}Q(s_0,t).
\end{align}
Consequently,
\begin{align}\notag
{II}&\ge C_1 z(s_0,t)Q(s_0,t)-C_2 s_0^{1/n}Q(s_0,t).
\end{align}

Clearly,
\begin{align*}
{III}&\ge -s_0[Q(s_0,t)+M]\\
&\ge -s_0^{\frac{1}{n}}Q(s_0,t)-s_0^{1-\frac2n}M.
\end{align*}

Combining the estimates for $I$, $II$, and $III$ yields
\begin{equation}\label{J ineq}
J_1'(t)\ge\left(C_1 z(s_0,t)-(C_2+1) s_0^{1/n}\right)Q(s_0,t)-C_3 s_0^{1-\frac2n},
\end{equation}
where $C_2=C_2(n,\chi,R,M,\mu)$ and
$C_3=\frac{2n(n-1)}m(1+M)+M$ are independent of $W$ and $s_0$.

For $t\in[0,1]$, we have $w(r,t)\ge We^{-t}\ge W/e$, so
\[
 Q(s_0,t)\ge\frac We z(s_0,t)\ge\frac{W}{es_0}J_1(t).
\]

By \eqref{J_1(0) low}, we may choose $s_0>0$ sufficiently small that
\begin{equation}\label{s_0 small}
s_0<\frac{C_1\mu W}{8e\Bunit}
\end{equation}
and
\[
\frac{C_1J_1(0)}{4s_0}>(C_2+1)s_0^{1/n},\qquad
\frac{C_1WJ_1(0)^2}{2es_0^2}>C_3s_0^{1-2/n}.
\]
These inequalities allow us to define
 \[
 T_2^*:=\sup\left\{T\le T_1^*\ \middle|\
 \begin{aligned}
 &\frac{C_1J_1(t)}{4s_0}>(C_2+1)s_0^{1/n},\\
 &\frac{C_1WJ_1(t)^2}{2es_0^2}>C_3s_0^{1-2/n}
 \quad(0\le t<T)
 \end{aligned}
 \right\}.
\]
For $t\in[0,T_2^*)$, \eqref{J ineq} yields
\begin{align*}
J_1'(t)&\ge\left(C_1 s_0^{-1}J_1(t)-(C_2+1) s_0^{1/n}\right)\frac{W}{es_0}J_1(t)-C_3 s_0^{1-\frac2n}\\
&\ge \left(C_1 (2s_0)^{-1}J_1(t)-(C_2+1) s_0^{1/n}\right)\frac{W}{es_0}J_1(t)+\left(\frac{C_1WJ_1(t)^2}{2es_0^2}-C_3 s_0^{1-\frac2n}\right)\\
&\ge \frac{C_1W}{4e s_0^2}J_1(t)^2> 0.
\end{align*}
Monotonicity and continuity yield $T_2^*=T_1^*=1$.

Integrating over $(0,1)$ gives
\[
 \frac{C_1W}{4es_0^2}
 \le J_1^{-1}(0)-J_1^{-1}(1)
 \le  J_1^{-1}(0).
\]
Together with \eqref{s_0 small}, this implies
\[
J_1(0)\le  \frac{4es_0^2}{C_1W}<\frac{4es_0}{C_1W}\cdot \frac{C_1\mu W}{8e\Bunit}= \frac{\mu s_0}{2\Bunit},
\]
contradicting \eqref{J_1(0) low}. Hence $\Tmax<\infty$.

Finally, if $u$ remained bounded on $(0,\Tmax)$, elliptic comparison
would give
$\|v(\cdot,t)\|_{L^\infty}\le
\|u(\cdot,t)\|_{L^\infty}\|w(\cdot,t)\|_{L^\infty}$.
The $w$ equation and Gr\"onwall's inequality would then bound $w$ on this
finite interval, contradicting the extensibility criterion in
Lemma~\ref{global_bounded_lemma}.  Therefore
\[
 \limsup_{t\nearrow\Tmax}\|u(\cdot,t)\|_{L^\infty}=\infty.
\]
\end{proof}

\subsection{Proof of Proposition~\ref{prop:high-blowup-restatement}: The case~$1<m\le2-2/n$}

A second difficulty in higher dimensions is that it is hard to construct a suitable blowing-up subsolution that effectively controls the diffusion term. We therefore turn to a weighted energy method. Since the earlier mass-based diffusion estimate does not directly cover $m>1$, we introduce a logarithmic energy functional.

The moment argument of Tao and Winkler \cite{TaoWinkler2021} for linear diffusion was adapted by Wang and Wang \cite{wang2025} to $0<m<1$. Extending this argument to $m>1$ requires a different estimate of the diffusion term
$$
n^2s^{2-\frac{2}{n}}(1+z_s)^{m-1}z_{ss}
$$
After multiplying by a suitable power of $s$ and integrating by parts, a weighted integral involving $(1+z_s)^m$ appears. For $0<m\le1$, suitable weighted inequalities control this integral using the mass bound; for $m>1$, the mass bound alone is insufficient.

To overcome this difficulty, rather than reducing the exponent $m$, we multiply the equation by a negative power of $z+1$ (for simplicity, we later choose $(z+1)^{-1}$). In this way, the integration by parts produces, together with the unfavorable negative term involving $(1+z_s)^m$, an additional positive term involving $(1+z_s)^{m+1}$, which helps control it. Combining this observation with Lemma~\ref{lem:w-weight decreasing}, we obtain the following estimate.

\begin{proposition}\label{n ge 4}
Assume that $n\ge 4$ and $1<m\le 2-\frac{2}{n}$. There are constants $\widetilde{C}_1=\widetilde{C}_1(\chi,R,M)>0$, $\widetilde{C}_2=\widetilde{C}_2(n,\chi,R,M)>0$, and $\widetilde{C}_3=\widetilde{C}_3(n,m,M)>0$, independent of $W$ and the concentration scale, such that
\begin{equation*}\left(\int_0^{s_0}\log(1+z)\dd s\right)_t\ge \left[\widetilde{C}_1 z(s_0,t)-(\widetilde{C}_2+1) s_0^{1/n}\right]Q(s_0,t)-(\widetilde{C}_3+1)s_0^{2-\frac{2}{n}-m}
\end{equation*}
for $0<s_0<\min\{1,R^n\}$ and $0<t<\Tmax$.
\end{proposition}
\begin{proof}

Multiplying both sides of \eqref{high-z-equation} by $(1+z)^{-1}$ and integrating over $(0,s_0)$, we obtain

\begin{align}
\left(\int_0^{s_0}\log(1+z)\dd s\right)_t
&\ge \int_0^{s_0}\frac{n^2s^{2-\frac2n}(1+z_s)^{m-1}z_{ss}}{1+z}\dd s\notag\\
&\quad+\chi\int_0^{s_0}\frac{z_s}{1+z}
\left(Q-n\int_0^{s^{1/n}}\rho^{n-1}v\dd\rho\right)\dd s\notag\\
&\quad-\int_0^{s_0}\frac{z+Q}{1+z}\dd s\notag\\
&=:\widetilde{I}+\widetilde{II}+\widetilde{III},
\label{high-log-decomposition}
\end{align}
where we have discarded the term involving $\varphi$.

Clearly,
\begin{equation}\label{est of III}
\widetilde{III}\ge -s_0[Q(s_0,t)+1].
\end{equation}

To estimate $\widetilde{II}$, we write it as

\begin{align*}
\widetilde{II}&=\chi\int_0^{s_0}\frac{z_s}{1+z}Q\dd s
-n\chi\int_0^{s_0}\frac{z_s}{1+z}
\left(\int_0^{s^{1/n}}\rho^{n-1}v\dd\rho\right)\dd s\\
&=:\widetilde{II}_1-\widetilde{II}_2.
\end{align*}

The following terms $II_1$ and $II_2$ arise from the proof for the case $0<m\le1$.

Direct computation gives
\begin{align}\notag
\widetilde{II}_1&\ge (1+M)^{-1} II_1\\\notag
&\ge \chi(1+M)^{-1}(1+e^R)^{-1}  z(s_0,t)Q(s_0,t)=:\widetilde{C}_1 z(s_0,t)Q(s_0,t).\label{II1 est}
\end{align}

We next estimate $\widetilde{II}_2$. As in the estimate of $II_2$, we obtain
\begin{align*}
 n\int_0^{s^{1/n}}\rho^{n-1}v(\rho,t)\dd\rho\le Cs^{1/n} Q(R^n,t),~~0<s<R^n,
\end{align*}
and hence, since $0\le z\le M$,
\begin{align}\notag
\widetilde{II}_2:&= \chi\int_0^{s_0}\frac{z_s}{1+z}
 \left(n\int_0^{s^{1/n}}\rho^{n-1}v\dd\rho\right)\dd s\\\notag
&\le C s_0^{\frac{1}{n}} Q(R^n,t)\int_0^{s_0}z_s\dd s\\\notag
 &\le C s_0^{1/n}Q(R^n,t)z(s_0,t).
\end{align}

For $\sigma_1\le s\le\sigma_2$, Lemma~\ref{lem:w-weight decreasing} yields
\begin{align*}
z_s(\sigma_1,t)Q_s(\sigma_2,t)&=z_s(\sigma_1,t)z_s(\sigma_2,t)w(\sigma_2^{\frac{1}{n}},t)\\
&\le z_s(\sigma_1,t)z_s(\sigma_2,t)e^{R}w(\sigma_1^{\frac{1}{n}},t)\\
&=e^R Q_s(\sigma_1,t)z_s(\sigma_2,t).
\end{align*}

Integrating both sides with respect to $\sigma_1$ over $(0,s_0)$ and with respect to $\sigma_2$ over $(s_0,R^n)$, we obtain
$$z(s_0,t)(Q(R^n,t)-Q(s_0,t))\le e^{R}Q(s_0,t)(z(R^n,t)-z(s_0,t))$$
which implies
\begin{equation*}
z(s_0,t)Q(R^n,t)\le e^R Q(s_0,t)z(R^n,t).
\end{equation*}

Therefore, we can then further estimate $\widetilde{II}_2$ as follows.
\begin{align}\notag
\widetilde{II}_2:&= \chi\int_0^{s_0}\frac{z_s}{1+z}
 \left(n\int_0^{s^{1/n}}\rho^{n-1}v\dd\rho\right)\dd s\\\notag
&\le C s_0^{\frac{1}{n}} Q(R^n,t)\int_0^{s_0}z_s\dd s\\\notag
 &\le C s_0^{1/n}Q(R^n,t)z(s_0,t)\\\notag
 &\le C e^Rs_0^{1/n}z(R^n,t)Q(s_0,t)\\
&\le Ce^R M s_0^{1/n}Q(s_0,t)=:\widetilde{C}_2 s_0^{1/n}Q(s_0,t).\label{II2 est}
\end{align}

Therefore, we have
\begin{align}\notag
\widetilde{II}&\ge \widetilde{C}_1 z(s_0,t)Q(s_0,t)-\widetilde{C}_2 s_0^{1/n}Q(s_0,t)
\end{align}

We then turn to the estimate of $\widetilde{I}$.  Integration by parts yields
\begin{align}\notag
\widetilde{I}&= \int_0^{s_0}\frac{n^2s^{2-\frac{2}{n}}\left[(1+z_s)^{m}\right]_{s}}{m(1+z)}\dd s\\\notag
 &=\frac{n^2s_0^{2-\frac{2}{n}}(1+z_s(s_0,t))^{m}}{m(1+z(s_0,t))}
 -\frac{n(2n-2)}{m}\int_0^{s_0}\frac{s^{1-\frac{2}{n}}(1+z_s)^{m}}{1+z}\dd s\notag\\
 &~~~+\frac{n^2}{m}\int_0^{s_0}\frac{s^{2-\frac{2}{n}}z_s(1+z_s)^{m}}{(1+z)^2}\dd s.
 \label{high-diffusion-identity}
\end{align}

We retain the positive integral in
\eqref{high-diffusion-identity} and obtain
\begin{align}\notag
 \widetilde{I}&\ge -\frac{2n(n-1)}{m}\int_0^{s_0}s^{1-\frac{2}{n}}(1+z_s)^{m}\dd s\notag\\\notag
 &~~~+\frac{n^2}{m}\int_0^{s_0}\frac{s^{2-\frac{2}{n}}(1+z_s)^{m+1}}{(1+z)^2}\dd s-\frac{n^2}{m}\int_0^{s_0}\frac{s^{2-\frac{2}{n}}(1+z_s)^{m}}{(1+z)^2}\dd s\\\notag
 &\ge -\frac{2n^2}{m}\int_0^{s_0}s^{1-\frac{2}{n}}(1+z_s)^{m}\dd s\notag\\\notag
 &~~~+\frac{n^2}{m(1+M)^2}\int_0^{s_0}s^{2-\frac{2}{n}}(1+z_s)^{m+1}\dd s-\frac{n^2}{m}\int_0^{s_0}s^{1-\frac{2}{n}}(1+z_s)^{m}\dd s\\
 &= -\frac{3n^2}{m}\int_0^{s_0}s^{1-\frac{2}{n}}(1+z_s)^{m}\dd s+\frac{n^2}{m(1+M)^2}\int_0^{s_0}s^{2-\frac{2}{n}}(1+z_s)^{m+1}\dd s
 \label{est I}
\end{align}
By H\"older's and Young's inequalities,
\begin{align}
\int_0^{s_0}s^{1-\frac2n}(1+z_s)^m\dd s
&=\int_0^{s_0}
\left[s^{2-\frac2n}(1+z_s)^{m+1}\right]^{\frac{m-1}{m}}
\left[s^{2-\frac2n-m}(1+z_s)\right]^{\frac1m}\dd s\notag\\
&\le\left(\int_0^{s_0}s^{2-\frac2n}(1+z_s)^{m+1}\dd s\right)^{\frac{m-1}{m}}
\left(\int_0^{s_0}s^{2-\frac2n-m}(1+z_s)\dd s\right)^{\frac1m}\notag\\
&\le
\left(\int_0^{s_0}s^{2-\frac2n}(1+z_s)^{m+1}\dd s\right)^{\frac{m-1}{m}}\cdot s_0^{\frac{2-\frac2n-m}{m}}
\left(\int_0^{s_0}(1+z_s)\dd s\right)^{\frac1m}\notag\\
&\le\varepsilon\int_0^{s_0}s^{2-\frac2n}(1+z_s)^{m+1}\dd s
+C_\varepsilon(1+M)s_0^{2-\frac2n-m},
\label{for I}
\end{align}
where $C_\varepsilon$ depends only on $m$ and $\varepsilon$, and we used $m\le2-2/n$ and
\begin{align*}
\int_0^{s_0}(1+z_s)\dd s&=s_0+z(s_0,t)\\
&\le 1+M.
\end{align*}
Choosing $\varepsilon=[3(1+M)^2]^{-1}$ in \eqref{for I} and using \eqref{est I} yields
\begin{equation}
\widetilde{I}\ge-\widetilde{C}_3s_0^{2-\frac2n-m},
\end{equation}
with $\widetilde{C}_3=\widetilde{C}_3(n,m,M)$ independent of $W$ and $s_0$.
Since $0\le2-2/n-m<1$, the same power also controls $s_0$.

Substituting the lower bounds for $\widetilde{I}$, $\widetilde{II}$, and $\widetilde{III}$ into \eqref{high-log-decomposition}, we obtain

\begin{align}\notag
\left( \int_0^{s_0}\log (1+z) ds\right)_t&\ge -\widetilde{C}_3s_0^{2-\frac{2}{n}-m}+ \widetilde{C}_1 z(s_0,t)Q(s_0,t)-\widetilde{C}_2 s_0^{1/n}Q(s_0,t)-s_0[Q(s_0,t)+1]\\\notag
&=\left[\widetilde{C}_1 z(s_0,t)-\widetilde{C}_2 s_0^{1/n}-s_0\right]Q(s_0,t)-\widetilde{C}_3s_0^{2-\frac{2}{n}-m}-s_0\\
&\ge \left[\widetilde{C}_1 z(s_0,t)-(\widetilde{C}_2+1) s_0^{1/n}\right]Q(s_0,t)-(\widetilde{C}_3+1)s_0^{2-\frac{2}{n}-m}.
\end{align}
This is exactly what we aim to prove.
\end{proof}

We now prove Proposition~\ref{prop:high-blowup-restatement}
for $1<m\le2-2/n$.

\begin{proof}[Proof of Proposition~\ref{prop:high-blowup-restatement} for $1<m\le2-2/n$]
Fix $\mu>0$ and $W$ as in the proposition. For
$0<s_0<\min\{1,R^n\}$ to be chosen below, take a smooth nonnegative radial
$u_0$ satisfying \eqref{data-ass}, supported in $B_{(s_0/2)^{1/n}}(0)$,
with $\int_\Om u_0=\mu$, and set $w_0\equiv W$.
Set
\begin{equation}\label{high-log-moment}
 J_2(t)=\int_0^{s_0}\log(1+z(s,t))\dd s.
\end{equation}
Suppose, to the contrary, that the solution is global. Then
\begin{equation}\label{upper bounds}
 0\le J_2(t)\le s_0\log(1+M),
\end{equation}
and
\[
 z(s_0,t)\ge\frac{J_2(t)}{s_0}.
\]
For $t\in[0,1]$, we have $w(r,t)\ge We^{-t}\ge W/e$, so
\begin{equation}\label{q lower bd}
 Q(s_0,t)=\int_0^{s_0}z_s(s,t)w(s^{1/n},t)\dd s
 \ge\frac We z(s_0,t)
 \ge\frac{W}{es_0}J_2(t).
\end{equation}

We next estimate the initial data. Since $u_0$ is supported in $B_{(s_0/2)^{1/n}}(0)$, we have
$$z(s,0)=n\int_0^{s^{1/n}}\rho^{n-1}u_0(\rho)\dd\rho=\frac{\mu}{\Bunit},~~ s\ge \frac{s_0}{2},$$ and hence
\begin{equation}\label{high-initial-moment}
J_2(0)\ge\int_{\frac{s_0}{2}}^{s_0}\log (1+z) ds=\frac {s_0}2\log\left(1+\frac{\mu}{\Bunit}\right)=:C_\mu s_0.
\end{equation}

For $t\in(0,1]$, Proposition~\ref{n ge 4} gives
\begin{equation}\label{new scalar2}
\left( \int_0^{s_0}\log (1+z) ds\right)_t\ge \left[\widetilde{C}_1 z(s_0,t)-(\widetilde{C}_2+1) s_0^{1/n}\right]Q(s_0,t)-(\widetilde{C}_3+1)s_0^{2-\frac{2}{n}-m}.
\end{equation}

Using \eqref{q lower bd}, we infer from \eqref{new scalar2} that
\begin{align}\notag
J_2'(t)&\ge \left[\widetilde{C}_1s_0^{-1}J_2(t)-(\widetilde{C}_2+1) s_0^{1/n}\right]Q(s_0,t)-(\widetilde{C}_3+1)s_0^{2-\frac{2}{n}-m}\\
&\ge\frac{\widetilde{C}_1W}{3e}\frac{J_2^2(t)}{s^2_0}+\left[\frac{\widetilde{C}_1J_2(t)}{3s_0}-(\widetilde{C}_2+1) s_0^{1/n}\right]Q(s_0,t)+\left[\frac{\widetilde{C}_1W}{3e}\frac{J_2^2(t)}{s^2_0}-(\widetilde{C}_3+1)s_0^{2-\frac{2}{n}-m}\right].
\label{new scalar3}
\end{align}

If $m=2-2/n$, the prescribed $W$ satisfies
\begin{equation}\label{high-W-choice}
W>\frac{3e(\widetilde{C}_3+1)}{\widetilde{C}_1C_\mu^2}.
\end{equation}
If $m<2-2/n$, keep any fixed $W>0$.
Since $2-2/n-m>0$ when $m<2-2/n$
and \eqref{high-W-choice} handles the endpoint, the independence of
$\widetilde{C}_1$ and $\widetilde{C}_3$ from $W$ allows $s_0$ to be chosen sufficiently small so that
\begin{gather*}
\frac{\widetilde{C}_1C_\mu}{3}>(\widetilde{C}_2+1)s_0^{1/n},\qquad
\frac{\widetilde{C}_1C_\mu^2W}{3e}>(\widetilde{C}_3+1)s_0^{2-\frac2n-m},\\
s_0<\frac{\widetilde{C}_1C_\mu W}{3e}.
\end{gather*}
Combining this with \eqref{high-initial-moment}, we further obtain
$$\frac{\widetilde{C}_1J_2(0)}{3s_0}-(\widetilde{C}_2+1) s_0^{1/n}>0$$
and
$$\frac{\widetilde{C}_1W}{3e}\cdot\frac{J_2^2(0)}{s^2_0}-(\widetilde{C}_3+1)s_0^{2-\frac{2}{n}-m}>0.$$

This ensures that the following $T^*$ is well defined.

\[
 T^*:=\sup\left\{T\le1\ \middle|\
 \begin{aligned}
 &\frac{\widetilde{C}_1J_2(t)}{3s_0}-(\widetilde{C}_2+1)s_0^{1/n}>0,\\
 &\frac{\widetilde{C}_1W}{3e}\cdot\frac{J_2^2(t)}{s_0^2}
   -(\widetilde{C}_3+1)s_0^{2-\frac2n-m}>0
   \quad(0\le t<T)
 \end{aligned}
 \right\}.
\]

We claim that $T^*=1$. Suppose, to the contrary, that $T^*<1$. Then, by the definition of $T^*$, at least one of
$$
\frac{\widetilde{C}_1J_2(T^*)}{3s_0}=(\widetilde{C}_2+1)s_0^{1/n}
$$
and
$$
\frac{\widetilde{C}_1W}{3e}\cdot\frac{J_2^2(T^*)}{s_0^2}
=(\widetilde{C}_3+1)s_0^{2-\frac{2}{n}-m}
$$
must hold.

On the other hand, by the definition of $T^*$ together with
\eqref{new scalar3}, we have
\begin{equation}\label{I last}
J_2'(t)\ge
\frac{\widetilde{C}_1W}{3e}\frac{J_2^2(t)}{s_0^2}\ge0.
\end{equation}
Hence $J_2$ is nondecreasing, and therefore
$$
\frac{\widetilde{C}_1J_2(T^*)}{3s_0}
\ge
\frac{\widetilde{C}_1J_2(0)}{3s_0}>(\widetilde{C}_2+1)s_0^{1/n},
$$
as well as
$$
\frac{\widetilde{C}_1W}{3e}\cdot\frac{J_2^2(T^*)}{s_0^2}
\ge
\frac{\widetilde{C}_1W}{3e}\cdot\frac{J_2^2(0)}{s_0^2}>
(\widetilde{C}_3+1)s_0^{2-\frac{2}{n}-m}.
$$
Thus both inequalities are strict, contradicting the fact that at least one of them must hold as an equality. Hence $T^*=1$.

We now know that \eqref{I last} holds for all $t\in[0,1]$, and hence we obtain
$$
-\frac{d}{dt}\left[J_2^{-1}(t)\right]\ge \frac{\widetilde{C}_1W}{3es_0^2}.$$

Integrating with respect to $t$ over $(0,1)$ and using \eqref{high-initial-moment}, we obtain
\begin{equation}
\frac{\widetilde{C}_1W}{3es_0^2}\le J_2^{-1}(0)-J_2^{-1}(1)\le \frac{1}{C_\mu s_0},
\end{equation}
which is impossible when $s_0<\frac{\widetilde{C}_1C_\mu W}{3e}$.
Thus $\Tmax<\infty$.

Finally, if $u$ remained bounded on $(0,\Tmax)$, elliptic comparison
would give
$\|v(\cdot,t)\|_{L^\infty}\le
\|u(\cdot,t)\|_{L^\infty}\|w(\cdot,t)\|_{L^\infty}$.
The $w$ equation and Gr\"onwall's inequality would then bound $w$ on this
finite interval, contradicting the extensibility criterion in
Lemma~\ref{global_bounded_lemma}.  Therefore
\[
 \limsup_{t\nearrow\Tmax}\|u(\cdot,t)\|_{L^\infty}=\infty.
\]

\end{proof}

\begin{remark}\label{rem:high-w0}
The argument also applies to radial $w_0$ satisfying \eqref{data-ass} and
$0<W_-\le w_0\le W_+$.  Indeed, write
\[
 w(t)=e^{t(\Delta_N-1)}w_0+\psi(t),\qquad
 \psi(t)=\int_0^t e^{(t-\sigma)(\Delta_N-1)}v(\sigma)\dd\sigma.
\]
The first term lies between $W_-e^{-t}$ and $W_+e^{-t}$, while
$(1-\Delta_N)\psi\ge0$.  The proof of
Lemma~\ref{lem:w-weight decreasing} therefore gives
\[
 w(r,t)\ge W_-e^{-t},\qquad
 w(r_2,t)\le
 \max\left\{e^R,\frac{W_+}{W_-}\right\}w(r_1,t)
 \qquad(r_1\le r_2).
\]
Replacing $e^R$ by this maximum in the comparison estimates and $W$ by $W_-$ in the lower bounds proves the same result. Here $C_1,C_2$ and $\widetilde{C}_1,\widetilde{C}_2$ also depend on $W_+/W_-$, while $C_3$ and $\widetilde{C}_3$ are unchanged. For $m<2-2/n$, any fixed such profile is admissible. For $m=2-2/n$, choose $W_-$ to satisfy \eqref{high-W-choice} with the resulting constants, keeping $W_+/W_-$ fixed.
\end{remark}

\section{Global boundedness for the fully parabolic system with \texorpdfstring{$n=3$}{n=3}}\label{sec:full}

Throughout this section, let $n=3$.

We first establish a more general time-weighted estimate for the interaction
term $uw$.

\begin{lemma}
\label{lem:weighted-interaction}
Let $\tau\in\{0,1\}$, and suppose that $(u,v,w)$ is a solution of
\eqref{sys}.  Then, for every $\delta>0$, there exists $C>0$ such that
\begin{equation*} \sup_{t\in(0,\Tmax)}
 \int_0^t e^{\delta(s-t)}\int_\Om u(x,s)w(x,s)\dd x\dd s\le C.
\end{equation*}
\end{lemma}

\begin{proof}
The case $\delta=1$ was established in
\cite{BaiLiLiuLuoPreprint}.  For self-containment, we include the argument
and then extend it to arbitrary $\delta>0$.
Multiplying the first equation in \eqref{sys} by $e^t$ and integrating
over $\Om\times(0,t)$, we obtain
\begin{align}\nonumber
 \int_0^t e^{s-t}\int_\Om u(x,s)w(x,s)\dd x\dd s
 &=e^{-t}\int_\Om u_0
  +\int_0^t e^{s-t}\int_\Om\varphi(x,s)\dd x\dd s
  -\int_\Om u(x,t)\dd x\\
 &\le \int_\Om u_0+|\Om|\varphi_*=:C^*.\label{delta=1}
\end{align}
Consequently, the asserted estimate holds whenever $\delta\ge 1$.

It remains to consider the case $\delta\in(0,1)$.  Choose an
integer $k\ge0$ such that $k<t\le k+1$.  Splitting the time integral and
using \eqref{delta=1}, we obtain
\begin{align*}
 \int_0^t e^{\delta(s-t)}\int_\Om uw\dd x\dd s
 &=\sum_{i=1}^{k}\int_{i-1}^i e^{\delta(s-t)}\int_\Om uw\dd x\dd s
   +\int_k^t e^{\delta(s-t)}\int_\Om uw\dd x\dd s\\
 &\le\sum_{i=1}^{k}e^{-\delta t}\int_{i-1}^i
 e^{s-(1-\delta)(i-1)}\int_\Om uw\dd x\dd s
 +e^{1-\delta}\int_k^t e^{s-t}\int_\Om uw\dd x\dd s\\
 &\le\sum_{i=1}^{k}e^{-\delta t}\cdot C^*e^i
 \cdot e^{-(1-\delta)(i-1)}+e^{1-\delta}C^*\\
 &=C^*e^{-\delta t}\cdot e
 \left(\frac{e^{\delta k}-1}{e^\delta-1}\right)+e^{1-\delta}C^*\\
 &\le \frac{C^*e}{e^\delta-1}+e^{1-\delta}C^*.
\end{align*}
Here, the second inequality follows from
$\int_0^i e^{s}\int_\Om uw\dd x\dd s\le C^* e^{i}$, which is itself a
consequence of \eqref{delta=1}.  When $0<t\le1$, we take
$k=0$, in which case the sum is empty.
\end{proof}

The preceding time-weighted estimate yields a uniform $L^\infty$ bound for
$w$.

\begin{lemma}
\label{lem:full-w-infinity}
The solution of \eqref{sys} with $\tau=1$ satisfies
\begin{equation*} \sup_{t\in(0,\Tmax)}\|w(\cdot,t)\|_{L^\infty(\Om)}<\infty.
\end{equation*}
\end{lemma}

\begin{proof}
Fix $T<\Tmax$ and consider $t\in[0,T]$.  The variation-of-constants
formulas for the second and third equations in \eqref{sys} read

\begin{equation}\label{full-v-duhamel}
 v(t)=e^{t(\Delta_N-1)}v_0
      +\int_0^t e^{(t-s)(\Delta_N-1)}
       \bigl(u(\cdot,s)w(\cdot,s)\bigr)\dd s
\end{equation}
and
\begin{equation}\label{full-w-duhamel}
 w(t)=e^{t(\Delta_N-1)}w_0
      +\int_0^t e^{(t-s)(\Delta_N-1)}v(s)\dd s.
\end{equation}
Substituting \eqref{full-v-duhamel} into \eqref{full-w-duhamel} and using the
semigroup property and Fubini's theorem, we obtain
\begin{align*}
 w(t)
 &=e^{t(\Delta_N-1)}w_0
   +\int_0^t e^{(t-s)(\Delta_N-1)}e^{s(\Delta_N-1)}v_0\dd s\\
 &\quad+\int_0^t\int_0^s
 e^{(t-s)(\Delta_N-1)}e^{(s-\sigma)(\Delta_N-1)}
 \bigl(u(\cdot,\sigma)w(\cdot,\sigma)\bigr)\dd\sigma\dd s\\
 &=e^{t(\Delta_N-1)}w_0+te^{t(\Delta_N-1)}v_0
   +\int_0^t\int_\sigma^t
    e^{(t-\sigma)(\Delta_N-1)}
    \bigl(u(\cdot,\sigma)w(\cdot,\sigma)\bigr)\dd s\dd\sigma\\
 &=e^{t(\Delta_N-1)}w_0+te^{t(\Delta_N-1)}v_0
      +\int_0^t(t-\sigma)e^{(t-\sigma)(\Delta_N-1)}
       \bigl(u(\cdot,\sigma)w(\cdot,\sigma)\bigr)\dd\sigma.
 \end{align*}
Lemma~\ref{lem:neumann-semigroup} yields
\begin{align}\nonumber
 \|w(\cdot,t)\|_{L^\infty}
 &\le e^{-t}\|w_0\|_{L^\infty}+te^{-t}\|v_0\|_{L^\infty}
      +\int_0^t(t-\sigma)\|e^{(t-\sigma)(\Delta_N-1)}
       \bigl(u(\cdot,\sigma)w(\cdot,\sigma)\bigr)\|_{L^\infty}\dd\sigma\\\nonumber
 &\le (1+t)e^{-t}\|w_0+v_0\|_{L^\infty}
       +\int_0^t(t-\sigma)e^{-(t-\sigma)}C(1+(t-\sigma)^{-\frac{3}{2}})\|
        uw(\cdot,\sigma)\|_{L^1}\dd\sigma\\
 &\le \|w_0+v_0\|_{L^\infty}
       +9C\int_0^t(t-\sigma)^{-\frac{1}{2}}e^{-\frac{2}{3}(t-\sigma)}\|
        uw(\cdot,\sigma)\|_{L^1}\dd\sigma,\label{w infty est}
\end{align}
where we used the elementary estimate
$s^{3/2}e^{-s/3}(1+s^{-3/2})\le 9$ for $s=t-\sigma>0$,
together with the nonnegativity of $w_0$ and $v_0$ and
$(1+t)e^{-t}\le1$.

For fixed $t>0$, let $\mu_t$ be the measure on $(0,t)$ defined by
$\dd\mu_t(\sigma)
:=e^{-\frac{2}{3}(t-\sigma)}\dd\sigma$.  In the calculation below, we
write $\dd\mu$ for $\dd\mu_t(\sigma)$ whenever no confusion can arise.
Applying Lemma~\ref{lem:weighted-interaction} with $\delta=\frac{2}{3}$, we obtain
\begin{align*}
 \int_0^t(t-\sigma)^{-\frac12}
 \|uw(\cdot,\sigma)\|_{L^1}\dd\mu
 &\le\left(\int_0^t(t-\sigma)^{-\frac34}
 \|uw(\cdot,\sigma)\|_{L^1}\dd\mu\right)^{\frac23}
 \left(\int_0^t\|uw(\cdot,\sigma)\|_{L^1}
 \dd\mu\right)^{\frac13}\\
 &\le C\left(\int_0^t(t-\sigma)^{-\frac34}
 \|uw(\cdot,\sigma)\|_{L^1}\dd\mu\right)^{\frac23}\\
 &\le C\left(\int_0^t(t-\sigma)^{-\frac34}
 \|u(\cdot,\sigma)\|_{L^1}
 \|w(\cdot,\sigma)\|_{L^\infty}\dd\mu\right)^{\frac23}\\
 &\le C\sup_{\sigma\in(0,t)}
 \|w(\cdot,\sigma)\|_{L^\infty}^{\frac23}.
\end{align*}
Here, the last inequality follows from Lemma~\ref{lem:mass} and the
estimate
$\int_0^\infty r^{-\frac34}e^{-\frac23r}\dd r<\infty$.
Substituting this estimate into \eqref{w infty est}, we further obtain
\begin{align*}
 \|w(\cdot,t)\|_{L^\infty}
 &\le \|w_0+v_0\|_{L^\infty}
       +C\sup_{\sigma\in(0,t)}\|w(\cdot,\sigma)\|_{L^\infty}^{\frac{2}{3}}\\
 &\le \|w_0+v_0\|_{L^\infty}
       +C\sup_{\sigma\in[0,T]}\|w(\cdot,\sigma)\|_{L^\infty}^{\frac{2}{3}}.
\end{align*}
Taking the supremum over $t\in[0,T]$, we obtain
\begin{equation*}
 \sup_{t\in[0,T]}\|w(\cdot,t)\|_{L^\infty}
 \le \|w_0+v_0\|_{L^\infty}
       +C\sup_{\sigma\in[0,T]}\|w(\cdot,\sigma)\|_{L^\infty}^{\frac{2}{3}},
\end{equation*}
which, by Young's inequality, yields a bound for
$\sup_{t\in[0,T]}\|w(\cdot,t)\|_{L^\infty}$ independent of $T$.
Since $T<\Tmax$ was arbitrary, this proves the lemma.
\end{proof}

The uniform signal bound now yields arbitrary finite $L^p$ bounds for $u$.
\begin{lemma}
\label{lem:full-u-Lp}
If $m>4/3$, then every solution of \eqref{sys} with $\tau=1$ satisfies
\begin{equation*} \sup_{t\in(0,\Tmax)}\|u(\cdot,t)\|_{L^p(\Om)}<\infty
 \qquad\text{for every }p>1.
\end{equation*}
\end{lemma}

\begin{proof}
Fix $p\ge2$.  Testing the first equation in \eqref{sys} with $pu^{p-1}$
and integrating by parts, we obtain
\begin{align}
 \frac{\dd}{\dd t}\int_\Om u^p
 &+p(p-1)\int_\Om
 u^{p-2}(1+u)^{m-1}|\nabla u|^2
 +p\int_\Om u^p+p\int_\Om u^pw\notag\\
 &=-\chi(p-1)\int_\Om u^p\Delta v
   +p\int_\Om\varphi u^{p-1}.
 \label{full-u-exact}
\end{align}

Generalized Gagliardo--Nirenberg inequalities allowing summability exponents
below one were established in \cite[Lemma~2.2]{Winkler2002MMAS}.  Applying
the bounded-domain form \cite[Lemma~3.1]{WinklerDjie2010} to
$u^{\frac{p+m-1}{2}}$ with lower-order exponent
$\frac{2}{p+m-1}$, and using Lemma~\ref{lem:mass}, we obtain
\begin{align}
 \|u\|_{L^{3(p+m-1)}}^{p+m-1}
 &=\left\|u^{\frac{p+m-1}{2}}\right\|_{L^6}^2 \notag\\
 &\le C_p\left\|\nabla u^{\frac{p+m-1}{2}}\right\|_{L^2}^2
 +C_p\left(
 \int_\Om
 \left(u^{\frac{p+m-1}{2}}\right)^{\frac{2}{p+m-1}}
 \right)^{p+m-1} \notag\\
 &=C_p\left\|\nabla u^{\frac{p+m-1}{2}}\right\|_{L^2}^2
 +C_p\left(\int_\Om u\right)^{p+m-1} \notag\\
 &\le C_p\left\|\nabla u^{\frac{p+m-1}{2}}\right\|_{L^2}^2+C_p \notag\\
 &\le C_p\int_\Om
 u^{p-2}(1+u)^{m-1}|\nabla u|^2+C_p.
 \label{full-coercivity}
\end{align}

Combining \eqref{full-u-exact} and \eqref{full-coercivity} and using
$w\ge0$, Young's inequality, and H\"older's inequality, we obtain
\begin{align*}
 &\frac{\dd}{\dd t}\int_\Om u^p+\int_\Om u^p
 +c_p\|u\|_{L^{3(p+m-1)}}^{p+m-1} \\
 &\le \frac{\dd}{\dd t}\int_\Om u^p+\int_\Om u^p
 +p(p-1)\int_\Om
 u^{p-2}(1+u)^{m-1}|\nabla u|^2+C_p \\
 &=-\chi(p-1)\int_\Om u^p\Delta v
   +p\int_\Om\varphi u^{p-1}
   -(p-1)\int_\Om u^p
   -p\int_\Om u^pw+C_p \\
 &\le \chi(p-1)\int_\Om u^p|\Delta v|
   +\left[(p-1)\int_\Om u^p+\int_\Om\varphi^p\right]
   -(p-1)\int_\Om u^p+C_p \\
 &\le \chi(p-1)\int_\Om u^p|\Delta v|
   +|\Om|\varphi_*^p+C_p \\
 &\le C_p\|u\|_{L^{2p}}^p\|\Delta v\|_{L^2}+C_p.
\end{align*}

Multiplication by $e^t$, integration over $(0,t)$, and division by $e^t$
give
\begin{align}
 \int_\Om u^p(\cdot,t)
 &+c_p\int_0^t e^{s-t}
 \|u(\cdot,s)\|_{L^{3(p+m-1)}}^{p+m-1}\dd s \notag\\
 &\le e^{-t}\int_\Om u_0^p
 +C_p\int_0^t e^{s-t}
 \|u(\cdot,s)\|_{L^{2p}}^p
 \|\Delta v(\cdot,s)\|_{L^2}\dd s
 +C_p\int_0^t e^{s-t}\dd s \notag\\
 &\le C_p+C_p\int_0^t e^{s-t}
 \|u(\cdot,s)\|_{L^{2p}}^p
 \|\Delta v(\cdot,s)\|_{L^2}\dd s.
 \label{full-u-energy}
\end{align}

Interpolation between $L^1(\Om)$ and $L^{3(p+m-1)}(\Om)$, together with
Lemma~\ref{lem:mass}, gives
\begin{align}
 \|u\|_{L^{2p}}^p
 &\le
 \|u\|_{L^{3(p+m-1)}}^{p\theta}
 \|u\|_{L^1}^{p(1-\theta)} \notag\\
 &\le C_p
 \left(
 \|u\|_{L^{3(p+m-1)}}^{p+m-1}
 \right)^{\frac{p\theta}{p+m-1}},
 \label{full-u-2p}
\end{align}
where $\theta\in(0,1)$ satisfies
\begin{equation}\label{alpha}
 \frac{1}{2p}
 =\frac{\theta}{3(p+m-1)}+1-\theta.
\end{equation}

By \eqref{alpha} and the assumption $m>\frac43$,
$\frac{p\theta}{p+m-1}=\frac{3(2p-1)}{2(3p+3m-4)}<1$.
Thus, \eqref{full-u-2p}, \eqref{full-u-energy}, and Young's inequality
yield
\begin{align*}
 \int_\Om u^p(\cdot,t)
 &+c_p\int_0^t e^{s-t}
 \|u(\cdot,s)\|_{L^{3(p+m-1)}}^{p+m-1}\dd s\\
 &\le C_p+C_p\int_0^t e^{s-t}
 \left(
 \|u(\cdot,s)\|_{L^{3(p+m-1)}}^{p+m-1}
 \right)^{\frac{p\theta}{p+m-1}}
 \|\Delta v(\cdot,s)\|_{L^2}\dd s\\
 &\le C_p+\frac{c_p}{2}\int_0^t e^{s-t}
 \|u(\cdot,s)\|_{L^{3(p+m-1)}}^{p+m-1}\dd s
 +C_p\int_0^t e^{s-t}
 \|\Delta v(\cdot,s)\|_{L^2}^{
 \frac{2(3p+3m-4)}{6m-5}}\dd s.
\end{align*}
Hence,
\begin{align}
 \int_\Om u^p(\cdot,t)
 &+\frac{c_p}{2}\int_0^t e^{s-t}
 \|u(\cdot,s)\|_{L^{3(p+m-1)}}^{p+m-1}\dd s\notag\\
 &\le C_p+C_p\int_0^t e^{s-t}
 \|\Delta v(\cdot,s)\|_{L^2}^{
 \frac{2(3p+3m-4)}{6m-5}}\dd s.
 \label{full-u-weighted}
\end{align}

Applying Lemma~\ref{lem:MR} to the $v$ equation and using the
$L^\infty$ bound for $w$ from Lemma~\ref{lem:full-w-infinity}, we obtain
\begin{align}
 \int_0^t e^{s-t}\|\Delta v\|_{L^2}^{
 \frac{2(3p+3m-4)}{6m-5}}\dd s
 &\le C+C\int_0^t e^{s-t}\|uw\|_{L^2}^{
 \frac{2(3p+3m-4)}{6m-5}}\dd s\notag\\
 &\le C+C\int_0^t e^{s-t}\|u\|_{L^2}^{
 \frac{2(3p+3m-4)}{6m-5}}\dd s\notag\\
 &\le C+C\int_0^t e^{s-t}
 \|u\|_{L^{3(p+m-1)}}^{
 \frac{2(3p+3m-4)\beta}{6m-5}}
 \|u\|_{L^1}^{
 \frac{2(3p+3m-4)(1-\beta)}{6m-5}}\dd s\notag\\
 &\le C+C\int_0^t e^{s-t}
 \left(
 \|u\|_{L^{3(p+m-1)}}^{p+m-1}
 \right)^{\frac3{6m-5}}\dd s,
 \label{full-v-MR}
\end{align}
where we used Lemma~\ref{lem:mass} and interpolation between $L^1(\Om)$
and $L^{3(p+m-1)}(\Om)$, with $\beta\in(0,1)$ satisfying
\begin{equation}\label{beta}
 \frac12=\frac{\beta}{3(p+m-1)}+1-\beta.
\end{equation}

Since $m>\frac43$, we have $\frac3{6m-5}<1$.  Hence, Young's inequality
applied to the last term in \eqref{full-v-MR} yields
\begin{equation}
 C_p\int_0^t e^{s-t}\|\Delta v\|_{L^2}^{
\frac{2(3p+3m-4)}{6m-5}}\dd s\le C+\frac{c_p}{4}\int_0^t e^{s-t}
 \|u\|_{L^{3(p+m-1)}}^{p+m-1}\dd s.
 \label{last-v-MR}
\end{equation}

Substituting \eqref{last-v-MR} into \eqref{full-u-weighted} and absorbing
gives the uniform $L^p$ bound for $p\ge2$.  Since $\Om$ is bounded,
H\"older's inequality extends the bound to all $p>1$.
\end{proof}

The preceding estimates yield global boundedness for the fully parabolic
system by a standard final iteration.

\begin{proof}[Proof of Theorem \ref{thm:main-parabolic}]
Lemma~\ref{lem:full-w-infinity} gives a uniform $L^\infty(\Om)$ bound for
$w$, while Lemma~\ref{lem:full-u-Lp} provides uniform $L^p(\Om)$ bounds
for $u$ for every finite $p>1$.  Thus $uw$ is bounded in every finite
$L^p(\Om)$, and the $v$ equation provides the regularity required by the
Moser--Alikakos iteration.  Applying
\cite[Lemma~A.1]{TaoWinkler2012JDE} and standard parabolic
regularity yields the asserted boundedness; see also
\cite[Lemma~4.3]{wang2025} for this standard concluding argument in the
chemotaxis--May--Nowak setting.  Finally, Lemma~\ref{global_bounded_lemma} implies
$\Tmax=\infty$.
\end{proof}

\end{document}